\documentclass[11pt]{article} 
\usepackage[utf8]{inputenc}
\usepackage[T1]{fontenc}
\usepackage[english]{babel}
\usepackage[paper=letterpaper, margin=1in]{geometry}
\usepackage{setspace}
\usepackage{amsmath, amssymb, amsthm, amsfonts, latexsym}
\usepackage{bbm}
\usepackage{mathtools}
\usepackage{algorithm}
\usepackage{algpseudocode}
\usepackage{authblk}

\usepackage{graphicx}
\usepackage{tikz}
\usetikzlibrary{calc, positioning, arrows.meta} 
\usepackage[dvipsnames]{xcolor}
\usepackage{xcolor}
\definecolor{satBlue}{HTML}{2173e8}
\definecolor{satGreen}{HTML}{0adc70}
\definecolor{satRed}{HTML}{aa5957}
\usepackage{subcaption} 
\usepackage{multirow}
\usepackage[shortlabels]{enumitem}
\usepackage{float}

\usepackage[numbers, sort&compress]{natbib}

\usepackage[colorlinks=true, linkcolor=blue, citecolor=red, urlcolor=blue]{hyperref}

\newtheorem{thm}{Theorem}[section]
\newtheorem*{thm*}{Theorem} 
\newtheorem{lem}[thm]{Lemma}
\newtheorem{prop}[thm]{Proposition}
\newtheorem{cor}[thm]{Corollary}

\newtheorem{defn}[thm]{Definition} 
\newtheorem{claim}[thm]{Claim}

\newtheorem*{recallthm}{Theorem \ref{\currentrecallref}}
\newenvironment{recall}[1]
  {\newcommand{\currentrecallref}{#1}\begin{recallthm}}
  {\end{recallthm}}

\newtheorem*{recallcor_inner}{Corollary \ref{\currentrecallcorref}}
\newenvironment{recallcor}[1]
  {\newcommand{\currentrecallcorref}{#1}\begin{recallcor_inner}}
  {\end{recallcor_inner}}

\newtheorem*{recalllem_inner}{Lemma \ref{\currentrecalllemref}}
\newenvironment{recalllem}[1]
  {\newcommand{\currentrecalllemref}{#1}\begin{recalllem_inner}}
  {\end{recalllem_inner}}

\newtheorem*{recallprop_inner}{Proposition \ref{\currentrecallpropref}}
\newenvironment{recallprop}[1]
  {\newcommand{\currentrecallpropref}{#1}\begin{recallprop_inner}}
  {\end{recallprop_inner}}

\algrenewcommand\algorithmicrequire{\textbf{Input:}}
\algrenewcommand\algorithmicensure{\textbf{Output:}}
\title{A Hypergraph Container Method for Spread SAT: Approximation and Speedup}
\author[1]{Zicheng Han}
\author[1]{Yupeng Lin}
\author[1,2]{Jie Ma}
\author[1,3]{Xiande Zhang}
\affil[1]{School of Mathematical Sciences, University of Science and Technology of China, Hefei, Anhui 230026, China.}
\affil[2]{Yau Mathematical Sciences Center, Tsinghua University, Beijing 100084, China.}
\affil[3]{Hefei National Laboratory, University of Science and Technology of China, Hefei 230088, China}
\date{\today}

\begin{document}

\maketitle

\begin{abstract}

We develop a hypergraph container method for the Boolean Satisfiability Problem (SAT) via the newly developed container results~[\href{https://link.springer.com/article/10.1007/s00493-026-00214-1}{Campos and Samotij (2026)}]. 
This provides an explicit connection between the extent of spread of clauses and the efficiency of container-based algorithms. 
Informally, the more evenly the clauses are distributed, the stronger the shrinking effect of the containers, 
which leads to faster algorithms for SAT.

To quantify the extent of spread, 
 we use a weighted point of view, in which a clause of size $s$ receives weight $p^s$ for some $0<p\le 1$.
In this way, we introduce the notion of $(\lambda,p)_k$-structure for SAT formulas, where $\lambda$ is the spread parameter and $k$ is the maximum size of clauses. 
By the almost-independence property of containers, we prove that for formulas with $(\lambda,p)_k$-structures,
one can distinguish between ``unsatisfiable formulas'' and ``formulas 
satisfying at least a $(1-\delta)$-fraction of clauses'' in sub-exponential time. 
This shows that sufficiently spread formulas are not worst-case instances for Gap-ETH.
Moreover, we show that the speedup is directly controlled by the spread parameter $\lambda$, 
yielding faster exact algorithms for SAT formulas containing a $(\lambda,p)_k$-structure. 
This result extends previous work~[\href{https://dl.acm.org/doi/pdf/10.1145/3564246.3585163}{Zamir (STOC 2023)}] to the non-uniform case.
\end{abstract}

\section{Introduction}
The Boolean Satisfiability Problem (SAT) asks whether a given Boolean formula admits a satisfying assignment. 
Typically, all formulas are in Conjunctive Normal Form (CNF)~\footnote{
A conjunction of clauses $\bigwedge_{i} K_i$, where each clause $K_i$ is a disjunction of literals $\bigvee_{j} l_{i,j}$. 
Throughout the paper, we assume that all CNF
formulas are subsumption-free.}, and when every clause has exactly $k$ literals, the problem is called $k$-SAT. 
SAT occupies a central place in complexity theory through its role in $\mathsf{NP}$-completeness, and the design of faster SAT algorithms has long been a major theme in algorithmic research.

From the viewpoint of exact exponential-time algorithms, the basic bound is the trivial $O^*(2^n)$~\footnote{Throughout this paper, we write $O^*(f(n)) := O(f(n)\cdot \mathrm{poly}(n))$.} exhaustive search. 
Let $\beta_k\in(1,2] $ be the smallest constant such that $k$-SAT can be solved in $(\beta_k+o(1))^n $ time. 
A main goal is to understand how far $\beta_k$ can be improved below 2. 
This question is closely related to the \textit{ Exponential Time Hypothesis} (ETH, see~\cite{impagliazzo2001complexity}) and the \textit{ Strong Exponential Time Hypothesis} (SETH, see~\cite{impagliazzo2001problems}), which provide a complexity-theoretic picture of worst-case SAT. 
Over the years, the PPSZ algorithm~\cite{paturi2005improved,hansen2019faster,hertli2014,scheder2021ppsz}(a randomized procedure based on variable ordering and logical resolution) has led to the best known upper bounds for $k$-SAT.

Rather than improving the worst-case running time for arbitrary SAT formulas, recent work of Zamir~\cite{zamir2023algorithmic} showed that a well-spread structure can be applied to obtain algorithmic speedups for $k$-SAT. 
Motivated by this perspective, we investigate a quantitative relationship between spread structure and tractability of SAT formulas, aiming to understand how structural properties turn into both approximation guarantees and speedups. 
In this line of research, the hypergraph container method plays an important role in exploiting the spread structure.

\subsection{The Hypergraph Container Method and SAT}
The container method is a powerful and widely-used tool in extremal combinatorics. Originally developed for graphs~\cite{kleitman1982number}, it was independently extended to hypergraphs by Balogh, Morris, and Samotij~\cite{balogh2015independent}, and by Saxton and Thomason~\cite{saxton2015hypergraph}. Since then, the method has been further refined and improved in a series of works~\cite{nenadov2024probabilistic,balogh2020efficient,saxton2016simple}. More recently, Campos and Samotij ~\cite{campos2024towards} obtained an efficient version with sharper quantitative bounds.
The core idea is that every independent set can be fully enclosed in one of a small number of vertex sets, called \textit{containers},  each of which induces only a few edges. Let $\mathcal{H}$ be a hypergraph.  
Briefly speaking, the classical container methods provide the resulting family of containers $\mathcal{C}$ satisfying:
\begin{itemize}
    \item[(1)] Every independent set of $\mathcal H$ is contained in some container
    $C \in \mathcal C$.
    \item[(2)] The number of containers is $2^{o(|V(\mathcal H)|)}$.
    \item[(3)] Each container has smaller size than $|V(\mathcal H)|$.
    \item[(4)] Each container is \textit{``almost independent''}: it contains
    only a small fraction of the edges of $\mathcal H$.
\end{itemize}

In theoretical computer science, the container method has found notable applications. 
In property testing~\cite{blais2024new, blais2025testing, seth2025tolerant}, it yields optimal sample complexity by analyzing independent sets and sparse subgraphs. 
In counting $k$-SAT functions~\cite{dong2022enumerating,balogh2023nearly}, it helps to show that nearly all $k$-SAT functions are unate. 
Recently, it has been used to accelerate algorithms for $k$-SAT~\cite{zamir2023algorithmic}.

The hypergraph container method applies to SAT via a natural correspondence: a formula $\varphi$ defines a hypergraph $\mathcal{H}_{\varphi}$ such that each satisfying assignment of $\varphi$ corresponds to an independent set in $\mathcal{H}_{\varphi}$.
Let $\varphi$ be a SAT formula over variables $X = \{x_1, \dots, x_n\}$. The associated hypergraph $\mathcal{H}_{\varphi}$ is defined as follows: the vertex set $V(\mathcal{H}_{\varphi}) := \{x_i, \neg x_i \mid 1 \le i \le n\}$ consists of all $2n$ literals, and the edge set $\mathcal{E}(\mathcal{H}_{\varphi}) := \{ \{\neg{\ell} \mid \ell \in K\} : K \text{ is a clause of } \varphi \}$. 
Figure~\ref{hypergraphfig} provides such an example.
In particular, each edge represents a combination that completely falsifies a clause. 
Avoiding it guarantees at least one true literal, thereby satisfying the clause. 
Thus, this construction induces a bijection between satisfying assignments $\alpha: X\to \{0,1\}$ and independent sets $I_{\alpha} :=\{x_i \mid \alpha(x_i) = 1\} \cup \{\neg x_i \mid \alpha(x_i) = 0\}$.

\begin{figure}[htbp]
\centering
\includegraphics[scale=0.15]{hypergraph_1.png}
\caption{The hypergraph $\mathcal{H}_\varphi$ induced by $\varphi = \textcolor{satBlue}{(\neg x_1\lor \neg x_2\lor \neg x_3)}\bigwedge\textcolor{satGreen}{(\neg x_1\lor x_2)}\bigwedge\textcolor{satRed}{(x_2\lor  x_3)}$.}
\label{hypergraphfig}
\end{figure}

To capture the edge distribution required by the container method, Zamir~\cite{zamir2023algorithmic} introduced the \textit{$(D,c,\epsilon)$-structure} for $k$-uniform hypergraphs, defined by degree conditions: (i) $|\mathcal{E}(\mathcal{H})| \geq D \cdot |V(\mathcal{H})|$, (ii) $\Delta_1(\mathcal{H}) \leq c\cdot D$, and (iii) $\Delta_2(\mathcal{H}) \leq c\cdot D^{1-\epsilon}$, where $\Delta_i(\mathcal{H})$ denotes the maximum number of edges containing a specific vertex set of size $i$. 
He showed that $k$-SAT formulas ``containing'' such a structure (i.e., there is a subset $\mathcal{E}' \subseteq \mathcal{E}(\mathcal{H}_{\varphi})$ such that $(V(\mathcal{H}_{\varphi}), \mathcal{E}')$ is a $(D, c, \epsilon)$-structure) allow for an algorithmic speedup. 

\begin{thm}[\cite{zamir2023algorithmic}, Theorem 5.13]
\label{Zamir}
    For every $k, \beta, c, \epsilon > 0$, there exist $D_0, \delta > 0$ such that the following hold for any $D > D_0$:
    If there exists an algorithm for $k$-SAT running in $O^*(\beta^n)$ time, then $k$-SAT formulas containing a $(D,c,\epsilon)$-structure can be solved in $O^* ((\beta-\delta)^n)$ time.
\end{thm}

Zamir first applied the container method to exploit how spread structure in $k$-SAT achieves algorithmic speedups. 
By bounding the edge density and codegrees via the $(D,c, \epsilon)$-structure, this approach ensures that the constraints are not overly concentrated on a small set of literals. 
Hence, it achieves a running time faster than the worst-case bound. 
While discussing the necessity of the $(D,c, \epsilon)$-structure, Zamir posed a series of open questions: Can the $(D, c, \epsilon)$-structure be weakened without losing the algorithmic speedup? 
If a $k$-SAT formula is allowed to contain clauses of length strictly less than $k$, would this benefit the solving algorithm? 
Furthermore, can sharper quantitative results be obtained?

When moving from $k$-SAT to general SAT formulas with clauses of varying sizes, the corresponding hypergraph is no longer uniform. 
To address Zamir's questions, our starting point is that such bounded-width SAT formulas should be modeled as \textit{$k$-bounded hypergraphs}, in which each edge has size at most $ k$. 
In fact, clauses of different lengths show different constraining effects. 
By introducing new structural concepts that distinguish clauses of different sizes, we redefine the notion of spread for these formulas and obtain sharper quantitative results.

\subsection{A Unified Weighted Container Method for SAT}
Our work is inspired by~\cite[Theorem A]{campos2024towards}. 
Although \cite[Theorem A]{campos2024towards} is stated for the uniform setting, the same algorithmic framework naturally extends to $k$-bounded hypergraphs. 
They assign each edge $E$ a weight $p^{|E|}$ to control edges of different sizes appearing during the container algorithm.
This approach aligns with a fundamental observation in SAT formulas: shorter clauses act as more rigid constraints than longer clauses. 
Accordingly, this weighting scheme is not merely a convenient tool, but a natural way to measure the influence of clauses with varying size.

Using this notion, we give a unified view of the container method for general SAT formulas. 
To characterize the extent of spread, we introduce the \textit{$(\lambda,p)_k$-structure} for non-uniform hypergraphs, which is motivated by the $(D,c,\epsilon)$-structure.
Given a hypergraph $\mathcal H$ and a number $p\in(0,1]$, each edge $E$ is assigned a weight $p^{|E|}$ and $w_p(\mathcal H)$ denotes the total weight of all edges. 
For $i\ge1$, the maximal $p$-codegree $\Delta_{i,p}(\mathcal H)$ is defined as the maximum total weight of edges containing a fixed set of size $i$, and the average $p$-weight $\bar w_p(\mathcal H)$ is defined as $w_p(\mathcal H)/|V(\mathcal H)|$.
We call $\mathcal{H}$ a $(\lambda,p)_k$-structure if $\mathcal{H}$ is $k$-bounded and satisfies: (i) $\bar{w}_p(\mathcal{H}) \geq p $, (ii) $\Delta_{1,p}(\mathcal{H}) \leq\lambda \bar{w}_p(\mathcal{H})$, and (iii) $\Delta_{2,p}(\mathcal{H}) \leq \lambda\bar{w}_p(\mathcal{H})p^k $. 
For fixed $p$, $\lambda$ serves as a spread parameter.
When $\mathcal{H}$ is uniform, the $(D, c, \epsilon)$-structure becomes a $(\lambda, p)_k$-structure by choosing appropriate parameters.
In particular, we call a SAT formula $\varphi$ a $(\lambda,p)_k$-structure if $\mathcal{H}_{\varphi}$ is a $(\lambda,p)_k$-structure. We later make this condition more concrete through the random \(k\)-SAT model as an example.

We establish the following weighted container method for SAT formulas with a $(\lambda,p)_k$-structure, which is the main result of this paper.
This leads to a quantitative characterization of an important insight: stronger spread not only yields more efficient algorithms but also characterizes tractable SAT formulas.

\begin{thm}\label{technical}
For each $k, \lambda>1$ and $0<\delta<1/4$, there exists $p'$ such that the following hold for each $0<p<p'$. 
If $\varphi$ is a $(\lambda,p)_k$-structure on $n$ variables, then we can find a collection $\mathcal{S}$ of subsets of
$V(\mathcal{H}_{\varphi})$ (called fingerprints), and a collection $\mathcal{C}$ of subsets of
$V(\mathcal{H}_{\varphi})$ (called containers), in $O^*(2^{2H(2k^2p/\delta)n})$ time with the following properties.
\begin{itemize} 
\item[(1)] If there exists a satisfying assignment $\alpha$, then there exists a container $C\in \mathcal{C}$ such that $I_{\alpha}\subseteq C.$
    \item[(2)] For each  fingerprint $S\in \mathcal{S}$, $|S|\leq \frac{4k^2pn}{\delta}$.  Consequently, $|\mathcal{C}| \leq 2^{2H(2k^2p/\delta)n}$.
    \item[(3)] Each $C\in \mathcal{C}$ satisfies $C\cap \{x_i,\neg x_i\}\neq \emptyset$ for each $i\in [n]$.
   \item[(4)] For each container $C\in \mathcal{C}$, $|C|
   \leq (2-1/\lambda)n$. Moreover, for the hypergraph induced by $\varphi$, \[w_p(\mathcal{H}_{\varphi}[C])\leq \delta \cdot w_p(\mathcal{H}_{\varphi}).\]
\end{itemize}
Here $H(x)=-x\log_2 x-(1-x)\log_2(1-x)$.
\end{thm}

Theorem~\ref{technical} captures a key observation: the container method can serve as a powerful divide-and-conquer strategy for SAT solving. 
Given a formula $\varphi$, we apply this method to decompose it into a family of reduced formulas $\{\varphi[C]\}_{C \in \mathcal{C}}$ (see Definition~\ref{def:partial_formula}). 
Thus, the original satisfiability question is then reduced to deciding whether at least one of these reduced formulas is satisfiable.
This decomposition process is both correct and useful: 
Theorem~\ref{technical} (1) guarantees the correctness of this decomposition: $\varphi \equiv \bigvee_{C \in \mathcal{C}} \varphi[C]$.
Theorem~\ref{technical}~(2) shows that the number of reduced formulas is well-controlled.
Theorem~\ref{technical} (3) ensures that each container $C$ may contain a satisfying assignment, since at least one literal is selected from every variable pair.
Theorem~\ref{technical} (4) describes some advantages of containers: the size of each container is upper bounded; and the total weight of edges in $C$ is controlled, which is called \textit{almost-independence}. 
The difference between Theorem~\ref{technical} and the typical container method is that we only collect containers that possibly include an assignment. 
Hence, $\mathcal{C}$ can be empty, which directly implies that $\varphi$ is unsatisfiable. 

\begin{figure}[htbp]
\centering
\includegraphics[scale=0.42]{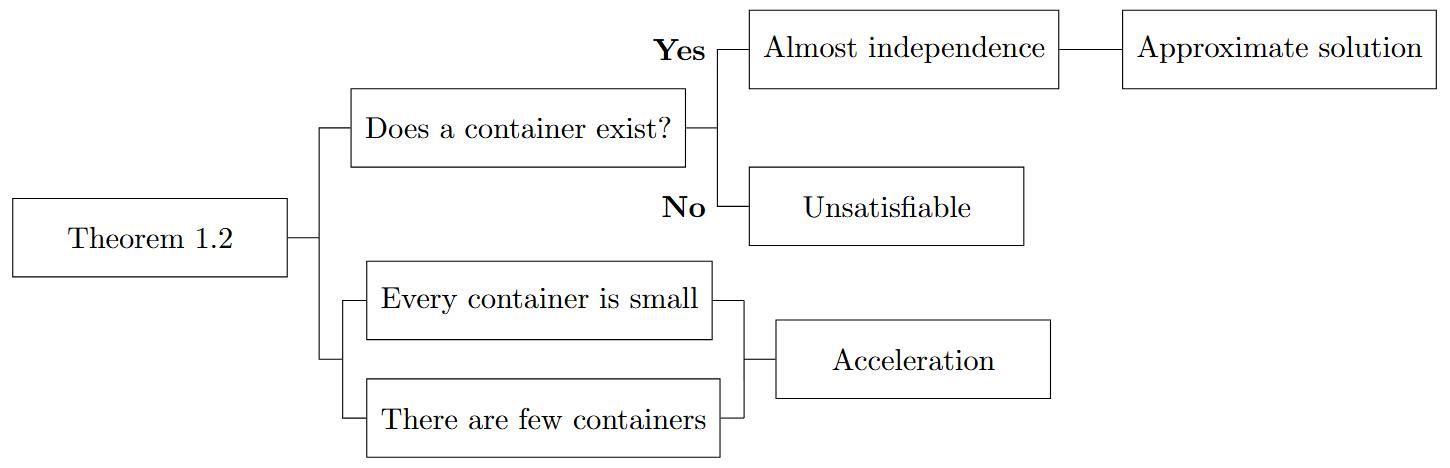}
\caption{How Theorem~\ref{technical} leads to our results. }
\label{way}
\end{figure}

These properties provide an effective problem-solving process. The number of containers is sub-exponential, and the size of each container is controlled by $\lambda$. In addition to taking less time, the algorithm also connects the speedup directly with the extent of spread. Moreover, the almost-independence property shows that assignments inside some container violate only a small fraction of clauses, and thus naturally leads to approximation algorithms. 
As a unified analytical tool, Theorem~\ref{technical} quantitatively establishes that the difficulty of a SAT formula is governed by how evenly its clauses are distributed: the more spread the formula is, the easier it becomes both to solve exactly and approximately. 
Figure~\ref{way} briefly shows how Theorem~\ref{technical} works.

\subsection{Approximate Satisfiability from Almost-Independence of Containers}

The Maximum Satisfiability problem (MAX-SAT) is the optimization variant of SAT: given a formula, find an assignment that maximizes the number of satisfied clauses. 
For MAX-E$k$-SAT, where each clause has exactly $k$ literals, the worst-case approximability is well understood:
H{\aa}stad~\cite{haastad2001some} proved that approximating MAX-E$k$-SAT within a factor of $(1 - 2^{-k} + \varepsilon)$ is NP-hard for any $\varepsilon > 0$.
Thus, in polynomial time, the random assignment threshold $(1 - 2^{-k})$ is optimal.
This naturally raises the question: can $(1 - 2^{-k})$ be broken if one allows algorithms to have super-polynomial time, such as sub-exponential time? 
For general formulas, the \textit{Gap Exponential Time Hypothesis} (Gap-ETH, see~\cite{impagliazzo2001complexity}) suggests that the answer remains negative: there exist constants $\delta, \varepsilon > 0$ such that no $O^*(2^{\varepsilon n})$-time algorithm can distinguish whether a 3-SAT formula is satisfiable or at most $(1-\delta)$-satisfiable~\footnote{Here, it means each assignment satisfies at most $(1-\delta)$-fraction of the clauses.}.

In this work, we introduce a new perspective, establishing a formal connection between the ``almost-independence'' property of hypergraph containers (Theorem~\ref{technical} (4)) and the solvability of MAX-E$k$-SAT. 
Our result shows that this worst-case barrier disappears under a structural spread assumption. 
In particular, for $k$-SAT formulas with $(\lambda,p)_k$-structures, we can either certify unsatisfiability or find an assignment satisfying at least $(1-\delta)$-fraction of clauses in sub-exponential time.

\begin{thm}\label{regularmax}
   For every $k,\lambda>1,0<\delta<1/4,$ and $\varepsilon>0$, there exists $p'$ such that for each $k$-SAT $(\lambda,p)_k$-structure $\varphi$ with $p<p'$, we can obtain one of the following outcomes in time $O^*(2^{\varepsilon n})$: 
       \begin{itemize}
    \item There is no satisfying assignment of $\varphi$, or
    \item There exists an assignment satisfying at least $(1-\delta)$-fraction of the clauses.
\end{itemize}
\end{thm}

Theorem~\ref{technical} (1) implies that if the family of containers $\mathcal{C}$ is empty, then $\varphi$ is unsatisfiable, since the independent set induced by an assignment must be contained in some container.
Otherwise, each container $C\in \mathcal{C}$ contains at least one literal from each variable pair, and thus contains an assignment. 
Theorem~\ref{technical} (4) implies that the ``almost-independence'' property upper-bounds the total weight of edges in $C$. 
Therefore, an assignment extracted from such a container violates only a small fraction of the clauses. 
This is how our result distinguishes between ``unsatisfiable'' and ``at least $(1-\delta)$-satisfiable'' for formulas with $(\lambda,p)_k$-structures.

We emphasize that this algorithmic result does not contradict known hardness assumptions such as Gap-ETH. 
While Gap-ETH assumes the $O^*(2^{\varepsilon n})$-hardness of distinguishing unsatisfiability from $(1-\delta)$-satisfiability for general worst-case formulas, our result applies solely to formulas with $(\lambda,p)_k$-structures. 
This confirms that formulas with sufficiently spread structure do not represent the hardest-case formulas predicted by Gap-ETH.

Analogous to~\cite{arora1995polynomial}, we also prove that it takes polynomial time to obtain one of the outcomes for formulas with a large clause-to-variable ratio.

\begin{cor}\label{densemax-k} 
For each $k>1$ and 
 $0<\delta<1/4,$ if the number of clauses in the $k$-SAT formula $\varphi$ is $\Theta(n^k)$,
 then we can obtain one of the following outcomes in polynomial time: 
       \begin{itemize}
   \item There is no satisfying assignment of $\varphi$, or
    \item There exists an assignment satisfying at least $(1-\delta)$-fraction of the clauses.
\end{itemize} 
\end{cor}

\subsection{Explicit Speedups for $(\lambda, p)_k$-Structures}
 Our weighted spread condition yields a direct consequence: formulas containing a \((\lambda,p)_k\)-structure can be solved faster than arbitrary SAT formulas, with an explicit dependence of the exponent on the spread parameter \(\lambda\).

\begin{thm}\label{dcepsilon}
  For each $ k, \lambda,\beta>1$ and $\varepsilon>0$, there exists $p'$ such that the following holds.
 If there exists an algorithm for SAT formulas running 
in $O^*(\beta^n)$ time, then SAT formulas containing a $(\lambda,p)_k$-structure with $p<p'$ can be solved in $O^*(\beta^{(1-1/\lambda+\varepsilon)n})$ time.
\end{thm}

This speedup stems directly from the structural properties of the containers established in Theorem~\ref{technical}. Specifically, combining the sub-exponential bound on the number of containers with the strict reduction of unassigned variables within them (Theorem~\ref{technical} (2)-(4)) yields an overall improvement in the running time exponent. 
The parameter \(\lambda\) measures how evenly the clauses are distributed, and our theorem shows that this structural information is reflected directly in the exponent through the factor (\(1-1/\lambda\)). 
Thus, the more evenly the clauses are distributed, the stronger the shrinking effect of the containers. 

Our theorem also recovers the $k$-SAT speedup obtained by Zamir~\cite{zamir2023algorithmic} as a special case. 
Zamir’s work already showed that container methods can accelerate \(k\)-SAT with a \((D,c,\epsilon)\)-structure. 
In our results, this appears as the uniform case of the more general weighted non-uniform viewpoint.
In particular, the $(D,c,\epsilon)$-structure is a special instance of the $(\lambda,p)_k$-structure under a suitable choice of parameters (see Corollary~\ref{cor:lambdap_dce}).

\begin{cor}\label{cor:lambdap_dce}
     For each $k,\beta,c>1$, $\varepsilon> 0$ and $0<\epsilon\leq 1$, there exists $D'$ such that the following holds.
If there exists an algorithm for $k$-SAT formulas running 
in $O^*(\beta^n)$ time, then a $k$-SAT formula containing a $(D,c,\epsilon)$-structure with $D>D'$ can be solved in $O^*(\beta^{(1-1/c+\varepsilon)n})$ time.
\end{cor}

The parameters $(1-1/\lambda)$ and $(1-1/c)$ suggest the same observation: the extent of speedup is tied to how evenly the edges are distributed in the graph. 
The application of optimal container methods~\cite{campos2024towards} allows us to make this dependence explicit.

Moreover, we demonstrate that a higher clause-to-variable ratio in SAT formulas yields a tighter upper bound on computational complexity. 

\begin{cor}\label{regular}
Let $p=1/(2n).$ 
For each $k, \beta >1$ and $d>0$, if there exists an algorithm for SAT formulas running in $O^*(\beta^n)$ time, then a SAT formula $\varphi$ containing a $k$-bounded subgraph $\mathcal{H}=(V(\mathcal{H}_{\varphi}),\mathcal{E}')$ with $w_p(\mathcal{H})\geq d$ can be solved in $O^*(\beta^{n-d/(2\Delta_{1,p}(\mathcal{H}))})$ time.
\end{cor}

\subsection{Outline of the paper}

Section~\ref{sec2} introduces preliminaries, including the SAT–hypergraph correspondence, decomposition of $\varphi$ and the definition of $p$-weights.
Section~\ref{sec3} derives the main results Theorems~\ref{regularmax} and~\ref{dcepsilon} from Theorem~\ref{technical}.
Section~\ref{sec4} gives a proof of Theorem~\ref{technical}, as well as Corollaries~\ref{densemax-k} and~\ref{regular}.
Section~\ref{sec5} contains concluding remarks.

\section{Preliminaries}\label{sec2}

\subsection{SAT Problems and Hypergraphs}
As discussed in the introduction, a SAT formula $\varphi$ naturally induces a hypergraph 
$\mathcal{H}_\varphi$ as follows. The vertex set consists of all literals 
$V(\mathcal{H}_\varphi)=\{x_i,\neg x_i:i\in[n]\}$, and for each clause $K$ we introduce an edge 
$E_K=\{\neg \ell:\ell\in K\}$.
This reduces a SAT problem to finding specific independent sets in its corresponding hypergraph. In fact,
satisfying assignments correspond exactly to independent sets that select exactly one literal from each variable pair. 
More precisely, each satisfying assignment $\alpha$ induces an independent set 
$I_\alpha=\{x_i\mid \alpha(x_i) = 1\} \cup \{\neg x_i \mid \alpha(x_i) = 0\}$, 
and $\alpha$ satisfies $\varphi$ if and only if $I_\alpha$ avoids every edge of $\mathcal{H}_\varphi$.
We denote by $\mathcal I^*(\mathcal H_{\varphi})$ the family of all such independent sets.

It is natural to reduce a SAT formula to a smaller SAT formula by focusing on some of its clauses. That is, for a SAT hypergraph $\mathcal{H}_\varphi$ induced by $\varphi$ and each subset $\mathcal{E}' \subseteq\mathcal{E}(\mathcal{H}_{\varphi})$, there exists a SAT formula $\phi$ whose induced SAT hypergraph $\mathcal{H}_\phi=(V(\mathcal{H}_\varphi),\mathcal{E}')$.

The following definition reduces a SAT formula to a smaller SAT formula by focusing on a designated subset of literals. Intuitively, variables outside $V'\subseteq V(\mathcal{H}_{\varphi})$
 are fixed (in particular, forced to take the value 0 whenever possible), while only the variables in $V'$ remain undetermined. Thus, any satisfying assignment on $V'$ can be viewed as a special satisfying assignment of the original SAT formula that is obtained by extending this partial assignment.

\begin{defn}\label{def:partial_formula}
For a vertex subset $V' \subseteq V(\mathcal{H}_{\varphi})$, the reduced formula $\varphi[V']$ denotes the formula $\varphi$ after partially assigning the following values to some of its variables.
For each variable $x_i$:
    \begin{itemize}
    \item[(1)] If $x_i \notin V'$ and $\neg x_i \in V'$, then $x_i$ is assigned to 0.
    \item[(2)] If $\neg x_i \notin V'$ and $x_i \in V'$, then $x_i$ is assigned to 1.
    \item[(3)] If both $x_i \in V'$ and $\neg x_i \in V'$, then $x_i$ remains unassigned.
    \item[(4)] If both $x_i \notin V'$ and $\neg x_i \notin V'$, then the formula $\varphi[V']$ is unsatisfiable.
    \end{itemize}
\end{defn}
The satisfiability of each reduced formula $\varphi[V']$ is closely related to the satisfiability of the original formula $\varphi.$ 
Moreover, following the rules above, we can eliminate the corresponding determined literals from further consideration, and thereby simplify the search for a satisfying assignment in $\varphi.$
The following lemmas formalize these observations in detail:
Lemma~\ref{lem:subgraph_sat} states the validity of the decomposition $\varphi \equiv \bigvee_{C \in \mathcal{C}} \varphi[C]$, while Lemma~\ref{lem:unassigned_bound} shows that solving each $\varphi[C]$ requires significantly less time than solving $\varphi$ itself. Here $\mathcal{C}$ is the collection of containers in Theorem~\ref{technical}.
\begin{lem}[\cite{zamir2023algorithmic}, Lemma 5.8] \label{lem:subgraph_sat}
    For each subset $V' \subseteq V(\mathcal{H}_{\varphi})$, if its reduced formula $\varphi[V']$ is satisfiable, then the original formula $\varphi$ is also satisfiable. Furthermore, if $\alpha$ is a satisfying assignment of $\varphi$ and its corresponding independent set $I_{\alpha} \in \mathcal{I}^*(\mathcal{H}_{\varphi})$ satisfies $I_{\alpha} \subseteq V' \subseteq V(\mathcal{H}_{\varphi})$, then $\varphi[V']$ is satisfiable.
\end{lem}

\begin{lem}[\cite{zamir2023algorithmic}, Lemma 5.9] \label{lem:unassigned_bound}
    Let $\gamma > 0$ and $\varphi$ be a SAT formula with $n$ variables. For any subset $V' \subseteq V(\mathcal{H}_{\varphi})$ with $|V'| \leq (1-\gamma)2n$, the reduced formula $\varphi[V']$ contains at most $(1-2\gamma)n$ unassigned variables.
\end{lem}

\subsection{The $p$-weight of Clauses}
In this subsection, we introduce some necessary notions to define the weight of clauses. 
Since our SAT formulas can contain clauses of different sizes, we assign different weights accordingly. 
Inspired by the $p$-weight notion introduced in~\cite{campos2024towards}, for a given real $p \in (0,1]$, we assign a weight of $p^s$ to a clause of size $s$. 
This weighting scheme reflects a simple intuition: longer clauses are easier to satisfy and should therefore carry less weight. 
Consequently, the notions of ``number of clauses'' and ``concentration of clauses'' must be reinterpreted under this view.
Recall that each SAT formula $\varphi$ can be represented as the SAT hypergraph $\mathcal{H}_{\varphi}$. The following definition formalizes this idea within the hypergraph context.

\begin{defn}
   Let $\mathcal{H}=(V(\mathcal{H}),\mathcal{E}(\mathcal{H}))$ be a hypergraph, and let $p \in (0,1]$ be a constant.
   \begin{itemize}
       \item The \textbf{$p$-weight of an edge $E$} is $p^{|E|}$. 
       \item For any subfamily $\mathcal{F} \subseteq \mathcal{E}(\mathcal{H})$, the \textbf{$p$-weight of $\mathcal{F}$} is $w_p(\mathcal{F}) = \sum_{E \in \mathcal{F}} p^{|E|}$. In particular, we abbreviate $w_p(\mathcal{E}(\mathcal{H}))$ as $w_p(\mathcal{H}) $.
       \item The \textbf{maximal $p$-codegree of $i$-sets} is defined as
       $\Delta_{i,p}(\mathcal{H}) := \max_{L \subseteq V, |L|=i} \sum_{E \supseteq L} p^{|E|}$.
       \item The \textbf{average $p$-weight} of $\mathcal{H}$ is defined as $\bar{w}_p(\mathcal{H})=w_p(\mathcal{H})/|V(\mathcal{H})|.$
   \end{itemize}
 
\end{defn}

The definitions $w_p(\mathcal{H})$, $\Delta_{i,p}$ and $\bar{w}_p$ serve as the weighted analogues of the edge number $|\mathcal{E}(\mathcal{H})|$, maximal codegree $\Delta_i$ and edge density $d$ in classical container methods.
For SAT formulas, the weight of a clause is defined as the weight of its corresponding edge. 

Finally, since no clause is contained in another, the hypergraph $\mathcal{H}_{\varphi}$ is an antichain. 
In this case, the LYM inequality~\cite{bollobas1965generalized,lubell1966short,yamamoto1954logarithmic,mesalkin1963generalization} provides a convenient bound to control the distribution of edges.

\begin{lem}[LYM inequality]\label{LYM}
    Let $\mathcal{F}\subseteq 2^{[n]}$ be an antichain, and let $a_s$ denote the number of sets of size $s$ in $\mathcal{F}$. Then
    \[
    \sum_{i=1}^n a_i\binom{n}{i}^{-1} \le 1.
    \]
\end{lem}

\section{From Theorem 1.2 to Other Main Results}\label{sec3}
In this section, we provide proofs of other main results derived from Theorem~\ref{technical}, whose proof is given in Section~\ref{sec4}.
Our results primarily focus on SAT problems whose clauses are relatively evenly distributed. 
We introduce the following definition to describe a generalized notion for evenly distributed non-uniform hypergraph.

\begin{defn}
\label{def:lambdap_structure}
    Let $\mathcal{H}$ be a hypergraph. 
    We say that $\mathcal{H}$ is a $(\lambda,p)_k$-structure if $\mathcal{H}$ is $k$-bounded and satisfies:
    \begin{itemize}
    \item [(1)] $\bar{w}_p(\mathcal{H})\geq p;$
       \item[(2)] $\Delta_{1,p}(\mathcal{H})\leq \lambda\cdot\bar{w}_p(\mathcal{H})$; and
       \item[(3)] $\Delta_{2,p}(\mathcal{H})\leq \lambda \cdot\bar{w}_p(\mathcal{H})\cdot p^{k}$.
       \end{itemize}
\end{defn}
Definition~\ref{def:lambdap_structure} ensures that in a $(\lambda,p)_k$-structure, no vertex or set has an excessively high weight.
It aligns with the common condition of the container method that edges cannot be overly concentrated.  The smaller $\lambda$ is, the more spread the hypergraph is. So we say $\lambda$ is a spread parameter. For SAT formulas,
we say $\varphi $ is a $(\lambda,p)_k$-structure if $\mathcal{H}_{\varphi}$ is a $(\lambda,p)_k$-structure,
and we say $\varphi$ contains a $(\lambda,p)_k$-structure if there exists a subset of edges $\mathcal{E}' \subseteq \mathcal{E}(\mathcal{H}_\varphi)$ such that the hypergraph $(V(\mathcal{H}_\varphi), \mathcal{E}')$ is a $(\lambda,p)_k$-structure.

 

\subsection{Approximate Solutions in Containers}\label{3.1}
The first result derived from Theorem~\ref{technical} is about the MAX-SAT problem. 
The almost-independent property of a container implies that the number of edges within it is extremely small compared to the original hypergraph, and this motivates the application of the container method to the MAX-SAT problem. 
In particular, if there exists a container $C$ such that $C\cap \{x_i,\neg x_i\}\neq \emptyset$ for each $i\in [n]$, then we can choose an assignment within $C$.
This chosen assignment is guaranteed to satisfy most of the clauses since $C$ itself is almost independent.

To formalize this intuition into the general MAX-SAT problem, we need to introduce a weighted version of MAX-SAT. 
This definition is well-motivated: shorter clauses are typically harder to satisfy and thus carry larger weight. 
Hence, the weighted version of MAX-SAT can be viewed as ``finding an assignment that minimizes the total weight of unsatisfied clauses.''
In this sense, the following lemma shows that when the weights are relatively evenly distributed, we can either find an approximate solution, or conclude that the problem is unsatisfiable.

\begin{lem}\label{maxsat}
   For each $k,\lambda>1,0<\delta<1/4$, there exists $p'$ such that for each $(\lambda,p)_k$-structure $\varphi$ with $p<p'$, we can obtain one of the following outcomes in time $O^*(2^{2H(2k^2p/\delta)n})$: 
       \begin{itemize}
    \item There is no satisfying assignment of $\varphi$, or
    \item There exists an assignment such that the weight of unsatisfied clauses is
    less than $\delta\cdot w_p(\mathcal{H}_{\varphi})$. 
\end{itemize}
\end{lem}

\begin{proof}[Proof of Lemma~\ref{maxsat} from Theorem~\ref{technical}]
   We apply Theorem~\ref{technical} with the given parameters to get a collection $\mathcal{C}$ of containers.
  If $\mathcal{C}=\emptyset$, then by Theorem~\ref{technical} (1), there is no satisfying assignment of $\varphi$. 
   Otherwise, $\mathcal{C}\neq\emptyset$.
For any such $C\in \mathcal{C}$, the induced subgraph $\mathcal{H}_{\varphi}[C]$ satisfies 
\[w_p(\mathcal{H}_{\varphi}[C])\leq \delta\cdot w_p(\mathcal{H}_{\varphi}).\]
By Theorem~\ref{technical} (3), $C$ contains at least one literal
from each variable pair $\{x_i,\neg x_i\}$.
Hence, for any assignment  $\alpha$ of $\varphi[C]$, its corresponding set (i.e., $ I_{\alpha}$, but it's not necessarily an independent set here) is contained in $C$. By the definition of $\varphi[C]$, each unsatisfied clause corresponds to an edge contained in $C$.
Consequently, the total weight of clauses falsified by $\alpha$ is at most $w_p(\mathcal{H}_{\varphi}[C])$, and hence is bounded by $\delta \cdot w_p(\mathcal{H}_{\varphi})$.
Finally, the total running time is $O^*(2^{2H(2k^2p/\delta)n})$ by Theorem~\ref{technical}. 
\end{proof}
When $\varphi$ is $k$-uniform, Lemma~\ref{maxsat} implies the conclusion with the MAX-E$k$-SAT problem.

\begin{recall}{regularmax}
   For every $k
,\lambda>1,0<\delta<1/4,$ and $\varepsilon>0$, there exists $p'$ such that for each $k$-SAT $(\lambda,p)_k$-structure $\varphi$ with $p<p'$, we can obtain one of the following outcomes in time $O^*(2^{\varepsilon n})$:
       \begin{itemize}
    \item There is no satisfying assignment of $\varphi$, or
    \item There exists an assignment  satisfying at least $(1-\delta)$-fraction of the clauses.
\end{itemize}
\end{recall}

\begin{proof}
By Lemma~\ref{maxsat}, for every $p<p'$, we can obtain one of the two outcomes in time $O^*(2^{2H(2k^2p/\delta)n})$.
Since the entropy function $H(x)\to 0$ as $x\to 0$, we may choose $p'$ sufficiently small such that $2H(2k^2p'/\delta)\leq \varepsilon$.
Then for every $0<p<p'$, it gives $2H(2k^2p/\delta)\leq 2H(2k^2p'/\delta)\leq \varepsilon$. 
Therefore, the time complexity is $O^*(2^{2H(2k^2p/\delta)n})\leq O^*(2^{\varepsilon n})$.
Since each edge has the same weight $p^k$, the outcome ``the weight of unsatisfied clauses in $\alpha$ is less than $\delta\cdot w_p(\mathcal{H}_{\varphi})$'' is equivalent to ``$\alpha$ satisfies at least $(1-\delta) $ fraction of clauses'', which completes our proof.
\end{proof}

The next results concern ``globally dense'' SAT formulas.
Here, globally dense refers to having a large clause-to-variable ratio,
and in the case of $k$-SAT, this corresponds to formulas containing
$\Theta(n^k)$ clauses.
Our corollaries show that global density is enough to distinguish ``unsatisfiable'' and ``existence of approximate solution'' in polynomial time.

\begin{cor}\label{densemax} 
Let $n$ be the number of variables and let $p=1/(2n).$ For each $k>1,d>0$ and 
 $0<\delta<1/4,$ if the $k$-bounded SAT formula $\varphi$ satisfies 
$w_p(\mathcal{H}_{\varphi})\geq d$, then we can obtain one of the following outcomes in polynomial time: 
       \begin{itemize}
   \item There is no satisfying assignment of $\varphi$, or
    \item There exists an assignment  such that the weight of unsatisfied clauses is
    less than $\delta\cdot w_p(\mathcal{H}_{\varphi})$.
\end{itemize} 
\end{cor}

The proof of Corollary~\ref{densemax} is postponed to Section 4.  From Corollary~\ref{densemax}, it is easy to get Corollary~\ref{densemax-k} as a special case.

\begin{recallcor}{densemax-k} 
For each $k>1$ and 
 $0<\delta<1/4,$ if the number of clauses in the $k$-SAT formula $\varphi$ is $\Theta(n^k)$,
 then we can obtain one of the following outcomes in polynomial time: 
       \begin{itemize}
   \item There is no satisfying assignment of $\varphi$, or
    \item There exists an assignment satisfying at least $(1-\delta)$-fraction of the clauses.
\end{itemize} 
\end{recallcor}

\subsection{Acceleration for Evenly Distributed SAT Problems}\label{3.2}
Our next result accelerates the search for a solution when the formula exhibits certain spread properties.
Based on Theorem~\ref{technical}, we establish the following theorem to achieve the speedup for $\varphi$ containing a $(\lambda,p)_k$-structure. 
 \begin{recall}{dcepsilon}
  For each $ k, \lambda,\beta>1$ and $\varepsilon>0$, there exists $p'$ such that the following holds.
 If there exists an algorithm for SAT formulas running 
in $O^*(\beta^n)$ time, then SAT formulas containing a $(\lambda,p)_k$-structure with $p<p'$ can be solved in $O^*(\beta^{(1-1/\lambda+\varepsilon)n})$ time.
\end{recall}

\begin{proof}[Proof of Theorem~\ref{dcepsilon} from Theorem~\ref{technical}:]
Let $\mathcal{H}=(V(\mathcal{H}_\varphi), \mathcal{E}')$ be the $(\lambda,p)_k$-structure contained in $\mathcal{H}_{\varphi}$, and let $\phi$ be the SAT formula with $\mathcal{H}_{\phi}=\mathcal{H}.$
We obtain $p''$ by applying Theorem~\ref{technical} to $\mathcal{H}_{\phi}$ with parameter $\delta = 1/5$, and let $p<p''$.
 This yields a family $\mathcal{C}$ of containers, and each of them has size less than $(2-1/\lambda)n$.
Since $\mathcal{H}_{\phi}$ is a subgraph of $\mathcal{H}_{\varphi}$, it implies that $\mathcal{I}^{*}(\mathcal{H}_{\varphi}) \subseteq \mathcal{I}^{*}(\mathcal{H}_{\phi})$.  Hence, for satisfying assignment $\alpha$ of $\varphi$, the independent set $I_\alpha\in \mathcal{I}^*(\mathcal{H}_{\varphi})\subseteq \mathcal{I}^*(\mathcal{H}_{\phi})$, and thus must be contained in some container of $\mathcal{H}_{\phi}$.

Therefore, to search for a satisfying assignment of $\varphi$, it suffices to search for a satisfying assignment of $\varphi[C]$ independently for each $C\in \mathcal{C}$ by Lemma~\ref{lem:subgraph_sat}.
It takes $O^*(2^{2H(10k^2p)n})$ time to identify all containers $C$.
The number of $\varphi[C]$s is upper-bounded by the number of containers, which is less than $2^{2H(10k^2p)n}.$
For each $\varphi[C]$, the number of unassigned variables is at most $(1 - 1/\lambda)n$ by Lemma~\ref{lem:unassigned_bound}. 
Using the assumed $O^*(\beta^n)$-time algorithm on $\varphi[C]$ takes time $O^*\big(\beta^{(1-1/\lambda)n}\big)$. 
Hence, the total running time is at most $O^*(2^{2H(10k^2p)n})+O^*(\beta^{(1-1/\lambda)n}\cdot 2^{2H(10k^2p)n})=O^*(\beta^{(1-1/\lambda)n}\cdot 2^{2H(10k^2p)n}).$

We choose $p'<p''$ such that 
$$2^{2H(10k^2p')n}\leq \beta^{\varepsilon n},$$
which is equivalent to 
\begin{equation}
2H(10k^2p')\leq \varepsilon\log_2\beta.
\label{thm1.5}
\end{equation}
Since $H(x)\to 0$ as $x\to 0$, we can choose $p'$ sufficiently small to satisfy \eqref{thm1.5}.
Hence for each $p\leq p'$, $H(10k^2p)\leq H(10k^2p')\leq \varepsilon \log_2 \beta/2$ always holds, completing the proof.
\end{proof}

The following corollaries demonstrate that ``faster problems'' in our results include a wide range of SAT problems, which guarantees the applicability of Theorem~\ref{dcepsilon}. 
The first example is the $(D,c,\epsilon)$-structure introduced by Zamir~\cite{zamir2023algorithmic}.
By showing that every $(D,c,\epsilon)$-structure is also a $(\lambda,p)_k$-structure, we demonstrate that the $(\lambda,p)_k$-structure captures a broad class of formulas.

\begin{lem}\label{lpdce}
    Each $k$-uniform $(D,c,\epsilon)$-structure with $D\geq 1,0<\epsilon\leq 1$ is a $(\lambda,p)_k$-structure for $\lambda = c$ and $ p=D^{-\epsilon/k}$.
\end{lem}
\begin{proof}
   Suppose that $\mathcal{H}$ is a $(D,c,\epsilon)$-structure. 
   By the definition, we have $D\leq \frac{ |\mathcal{E}(\mathcal{H})|}{|V(\mathcal{H})|}$,
   and the hypergraph $\mathcal{H}$ satisfies
\[\Delta_1\leq \frac{c\cdot |\mathcal{E}(\mathcal{H})|}{|V(\mathcal{H})|}\quad \mbox{and}\quad\Delta_2\leq c\left(\frac{ |\mathcal{E}(\mathcal{H})|}{|V(\mathcal{H})|}\right)D^{-\epsilon}.\]
Now we choose $p=D^{-\epsilon/k}$. 
Multiplying both sides of the degree bounds by $p^k$ yields the desired inequalities for 
$\Delta_{1,p}(\mathcal{H})$ and $\Delta_{2,p}(\mathcal{H})$ in the definition of $(\lambda,p)_k$-structure. 
Moreover, $$w_p(\mathcal{H})=|\mathcal{E}(\mathcal{H})|\cdot p^k=|\mathcal{E}(\mathcal{H})|D^{-\epsilon}\geq D^{1-\epsilon}|V(\mathcal{H})|.$$
By the condition $D\geq 1, $ we have $D^{1-\epsilon}\geq D^{-\epsilon/k}$, and it implies $\bar{w}_p(\mathcal{H})\geq p.$
Consequently, this shows that  $\mathcal{H}$ is a $(c,D^{-\epsilon/k})_k$-structure.   
\end{proof}
By Lemma~\ref{lpdce}, we improve Zamir's result \cite{zamir2023algorithmic} as a special case in Corollary~\ref{cor:lambdap_dce}, which explicitly characterizes the speedup by the spread parameters.
\begin{recallcor}{cor:lambdap_dce}
For each $k,\beta,c>1$, $\varepsilon> 0$ and $0<\epsilon\leq 1$, there exists $D'$ such that the following holds.
If there exists an algorithm for $k$-SAT formulas running 
in $O^*(\beta^n)$ time, then a $k$-SAT formula containing a $(D,c,\epsilon)$-structure with $D>D'$ can be solved in $O^*(\beta^{(1-1/c+\varepsilon)n})$ time.
\end{recallcor}

\begin{proof} We apply Theorem~\ref{dcepsilon} to  $k,\beta, \varepsilon$ and $\lambda=c$ to get the number $p'$. 
Let $\mathcal{H}$ be the $(D,c,\epsilon)$-structure in $\mathcal{H}_{\varphi}$.
By Lemma~\ref{lpdce}, it is a $(\lambda,p)_k$-structure for $\lambda=c,p=D^{-\epsilon/k}$.
We choose $D'$ large enough so that $(D')^{-\epsilon/k} \leq p'$. 
Then for any $D \geq D'$,
$p= D^{-\epsilon/k}\leq (D')^{-\epsilon/k}\leq p'.$
Thus $\mathcal{H}$ is a $(\lambda,p)_k$-structure with $\lambda=c, p =D^{-\epsilon/k}\leq p'$, and all conditions of Theorem~\ref{dcepsilon} are satisfied. 
\end{proof}
 
The second example also requires SAT formulas to be globally dense.   
Similarly, we show a better time complexity for $\mathcal{H}_{\varphi}$ in which the clause-to-variable ratio is a constant. 
The condition here is slightly different from the $(\lambda,p)_k$-structure.

\begin{recallcor}{regular}
Let $p=1/(2n).$ 
For each $k, \beta >1$ and $d>0$, if there exists an algorithm for SAT formulas running in $O^*(\beta^n)$ time, then a SAT formula $\varphi$ containing a $k$-bounded subgraph $\mathcal{H}=(V(\mathcal{H_{\varphi}}),\mathcal{E}')$ with  $w_p(\mathcal{H})\geq d$ can be solved in $O^*(\beta^{n-d/(2\Delta_{1,p}(\mathcal{H}))})$ time.
\end{recallcor}
The proof of  Corollary~\ref{regular} is postponed to Section 4.

In particular, the random \(k\)-SAT model provides a useful sanity check for the scope of \((\lambda,p)_k\)-structures. 
The model \(\Phi=\Phi_k(n,m)\) is defined by
choosing \(m\) distinct clauses uniformly at random from all 
\(k\)-clauses on variables \(x_1,\ldots,x_n\). 
For the associated SAT hypergraph $\mathcal{H}_{\Phi}$, the conditions of \((\lambda,p)_k\)-structure are
reduced to (set \(d=\frac{m}{2n}\)): 
\begin{equation}
    d\ge p^{1-k},\qquad
\Delta_1(\mathcal{H}_{\Phi})\le \lambda d,\qquad
\Delta_2(\mathcal{H}_{\Phi})\le \lambda d p^k,
\label{random_k_sat}
\end{equation}
where \(\Delta_1\) and \(\Delta_2\) denote the maximal degree and
maximal codegree.
The following proposition shows that, for fixed \(k,p,\lambda\), the scale \(m\asymp n\log n\) is the natural threshold for this random model to satisfy the \((\lambda,p)_k\)-condition.
\begin{prop}\label{bigrandom}
For each \(k\ge3\), \(p\in(0,1]\), and \(\lambda>k\):
\begin{itemize}
    \item[(1)] There exists a constant
\(C_0=C_0(k,p,\lambda)>0\) such that if \( d\ge C_0\log n,\)
then \(\mathcal H_\Phi\) is a \((\lambda,p)_k\)-structure w.h.p.
as \(n\to\infty\).
\item[(2)] If  $d= o(\log n)$, then $\mathcal{H}_{\Phi}$ is not a
\((\lambda,p)_k\)-structure w.h.p. as \(n\to\infty\). 
\end{itemize}  
\end{prop}
The proof of Proposition~\ref{bigrandom} is shown in the Appendix.

\section{Proof of Theorem~\ref{technical}}\label{sec4}
In this section, we present the container algorithm and further apply it to prove Theorem~\ref{technical}, along with Corollary~\ref{regular} and Corollary~\ref{densemax}. In Algorithm~\ref{alg:hypergraph_container}, we need the following notation.
Given a hypergraph $\mathcal H$, let $\mathcal H^s$ denote the family of edges in $\mathcal H$ with size $s$, and let $\mathcal H^{>s}$ (resp.\ $\mathcal H^{<s}$) denote the family of
edges whose sizes are strictly greater than $s$ (resp.\ strictly less than $s$).
For each subset $L \subseteq V(\mathcal H)$, define
$\partial_L \mathcal H=\{E\setminus L \mid L\subseteq E\in \mathcal E(\mathcal H)\}$.
These notations help to distinguish edges by size and describe the concentration of edges on a given vertex set.

Algorithm~\ref{alg:hypergraph_container} is adapted from the container method introduced in  \cite[Subsection~3.1]{campos2024towards}, which outputs the fingerprint and container for an independent set of a hypergraph. Here we modify some steps to suit SAT problems better. Although the input $\mathcal{H}$ is assumed to be $k$-uniform in~\cite{campos2024towards}, the same algorithm can be applied to $k$-bounded $\mathcal{H}$ in a natural way.
Moreover, let $\langle \mathcal{H}\rangle=\{F\mid \text{ there exists some } E\in \mathcal E(\mathcal H) \text{ such that } E\subseteq F\}$.
The stopping condition is changed to coincide with the remaining weights for MAX-SAT, and steps 16-19 are added to distinguish whether the resulting container is possible to contain an assignment.

\begin{algorithm}[htbp]
\caption{Hypergraph Container Algorithm}
\label{alg:hypergraph_container}
\begin{algorithmic}[1]
\Require $k$-bounded hypergraph $\mathcal{H}$, independent set $I\subseteq V(\mathcal{H})$~\footnotemark, constant $0<\delta<1/4$, and weight $0<p<\delta/(4k^2)$.
\Ensure Fingerprint $S$, container $C$, and hypergraph $\mathcal{G}$.

\State Let $\mathcal{H}_0 = \mathcal{H}$, $S_0 = \emptyset$.
\For{$i = 0, 1, \ldots$}
    \State Let $C_i = \{v \in V \mid \{v\} \notin \mathcal{H}_i\}$.
    \If{$w_p(\mathcal{H}_i^{>1}) \leq \delta p|C_i|$}
        \State Set $J = i$ and \textbf{STOP}.
    \EndIf
    \State Find the smallest integer $s\in \{2,\dots,k\}$ for which there exists a non-empty set $L \subseteq V(\mathcal{H})$ such that $L \notin \langle \mathcal{H}_i\rangle$ and $w_p(\partial_L (\mathcal{H}_i^s)) > \delta /k$. Denote the chosen integer by $s_i.$
    \State Among all such sets $L$ for $s_i$, choose one that is maximal under inclusion; if several maximal sets exist, pick the one with the smallest lexicographic order. Denote the chosen set by $L_i$.
    \If{$L_i \subseteq I$}
        \State Let $\mathcal{F}_i = \partial_{L_i}\mathcal{H}_i^{s_i}$ and $S_{i+1} = S_i \cup L_i$.
    \Else
        \State Let $\mathcal{F}_i = \{L_i\}$ and $S_{i+1} = S_i$.
    \EndIf
    \State Let $\mathcal{H}_{i+1} = (\mathcal{H}_{i} \setminus \langle \mathcal{F}_i \rangle) \cup \mathcal{F}_i$.
\EndFor
\If{there exists $x \in V$ such that both $x$ and $\neg x$ are not in $C_J$}
    \State Set $S_J, C_J, \mathcal{H}_J= \text{\textsc{Error}}$.
\Else
    \State $S_J,C_J,\mathcal{H}_J$ remain unchanged.
\EndIf
\State \Return $S = S_J$, $C = C_J$, $\mathcal{G} = \mathcal{H}_J^{>1}$.
\end{algorithmic}
\end{algorithm}
\footnotetext{By the fingerprint property (Lemma~\ref{3.11}), every independent set is contained in some container determined by a fingerprint of size at most $ 4k^2pn/\delta$.
Hence, it suffices to enumerate all independent sets of size at most this bound.}

Recall that the first condition of Definition~\ref{def:lambdap_structure} is \(\bar w_p(\mathcal H)\ge p\). 
It matches the stopping condition of Algorithm~\ref{alg:hypergraph_container}, and should be viewed as a minimal non-vacuity condition rather than a density assumption.
Steps 2-15 in Algorithm~\ref{alg:hypergraph_container} iteratively search for a set $L_i$ of high weight until the stopping condition holds. 
According to whether $L_i$ is contained in $I$, the algorithm records $L_i$ as part of the fingerprint, or turns $L_i$ itself into a new edge.
The iteration rule for $\mathcal{H}_i$ always maintains it as an antichain and pushes the remaining structure toward the stopping condition.  Steps 16-19 aim to filter out all bad containers $C_J$: there exists $x_i$ such that $C_J\cap \{x_i,\neg x_i\}=\emptyset$. 
Such a container cannot encode any satisfying assignment, since an assignment must include either $x_i$ or $\neg x_i$ for each variable $x_i$.
Consequently, if $C_J$ is eliminated for every input, then we conclude that the SAT formula $\varphi $ is unsatisfiable.
If this does not happen, it remains to estimate the weight of the edges contained in $\mathcal H[C]$.
The output hypergraph $\mathcal G$ helps this estimation.

Given Algorithm~\ref{alg:hypergraph_container}, we next show some key properties to check its efficiency. 
We begin by verifying that the algorithm is well-defined and terminates in finite time. The proofs for this non-uniform case are almost the same as those for the uniform case in~\cite{campos2024towards}, 
so we move them to  Appendix~\ref{Appendix} and give a brief explanation below.

\begin{lem}[\cite{campos2024towards}, Lemma 3.1-3.3]\label{3.1,3.2,3.3}
For each iteration step $i$ in Algorithm~\ref{alg:hypergraph_container}, the following properties hold:
\begin{itemize}
    \item There always exist desirable $s_i$ and $L_i$.
    \item $\mathcal{H}_i$ is an antichain. 
    \item $\langle\mathcal{H}_i\rangle\subsetneq \langle\mathcal{H}_{i+1}\rangle$, and $\mathcal{F}_i\cap\langle\mathcal{H}_{i}\rangle=\emptyset$.
\end{itemize}  
\end{lem}

Lemma~\ref{3.1,3.2,3.3} ensures that Algorithm~\ref{alg:hypergraph_container} is well-defined. The existence of $s_i$ and $L_i$ guarantees that the algorithm runs to the next step. 
Moreover, since each $\langle\mathcal{H}_i\rangle\subsetneq \langle\mathcal{H}_{i+1}\rangle$, the total number of steps cannot exceed $2^{2n}=4^n.$
This shows that the algorithm will finally terminate. 

A key point of the container algorithm is: to find all containers, it suffices to enumerate all small independent sets.
This is achieved via steps 9-12 by generating a \textit{fingerprint}, which is a small subset of the independent set such that any two independent sets sharing the same fingerprint lie in the same container.
Consequently, a small number of inputs implies that we spend less time generating all containers. These properties are summarized in the following 
 lemma, whose proof  is also moved to Appendix~\ref{Appendix}.

\begin{lem}[\cite{campos2024towards}, Lemma 3.5, 3.11]\label{3.11}
   For each independent set $I\in \mathcal{I}(\mathcal{H})$, $S_J\subseteq I\subseteq C_J$, where $J$ is the stopping index in Algorithm~\ref{alg:hypergraph_container}.
Moreover, for two different independent sets $I$ and $I'$ of $\mathcal{H}_0$, if $S=S'$, then $C=C'$ and $\mathcal{G}=\mathcal{G}'$.
\end{lem}

Next, we bound the size of each fingerprint via evaluating the concentration of short edges. In the original algorithm of Campos and Samotij, all edges of size less than \(k\) in $\mathcal H_i$ are newly generated, so one can control them by directly bounding \(w_p(\partial_L\mathcal H_i^s)\). 
Since here $\mathcal H_0$ may already contain edges with size less than $k$, we instead bound the partial weight of newly generated edges by introducing 
\[\mathcal{D}_i := \mathcal{H}_i \setminus \mathcal{H}_0,\] i.e., the collection of edges that have been added to the hypergraph $\mathcal{H}_i$ excluding those originally contained in $\mathcal{H}_0$. 
The next lemma provides a bound of \(w_p(\partial_L\mathcal D_i^s)\), indicating that the algorithm will not cause weight concentration of newly generated short edges.

\begin{lem}\label{upper} For each $i$, and each subset $L\notin \langle \mathcal{H}_i\rangle$ and $s$ satisfying
  $1\leq |L|<s\leq k-1$, $$w_p(\partial_{L}\mathcal{D}_i^{s})\leq \frac{2\delta}{k}.$$
\end{lem}
\begin{proof}
   We prove this lemma by induction on $i$. 
   For the base case $i=0,$ we have $\mathcal{D}_0=\mathcal{H}_0\backslash\mathcal{H}_0=\emptyset$ and these inequalities naturally hold.
   Suppose that these inequalities hold for integer $i$.
    The fact $\mathcal{H}_{i+1}\subseteq \mathcal{H}_{i}\cup \mathcal{F}_i$
     implies that $\mathcal{D}_{i+1}\subseteq \mathcal{D}_{i}\cup \mathcal{F}_i$. Hence for each $L$ and $s>|L|$,
     \[w_p(\partial_{L}\mathcal{D}_{i+1}^{s})\leq w_p(\partial_{L}\mathcal{D}_i^{s})+w_p(\partial_L \mathcal{F}_{i}^{s}).\]
     If $\partial_L \mathcal{F}_{i}^{s}=\emptyset,$ then $w_p(\partial_{L}\mathcal{D}_{i+1}^{s})\leq w_p(\partial_{L}\mathcal{D}_i^{s})\leq 2\delta /k$. 
     Otherwise, we have $\mathcal{F}_{i}^{ s}\neq \emptyset$. 
     Since each edge in $\mathcal{F}_i$ must be a proper subset of some edge in $\mathcal{H}_i^{s_i}$, it implies that $s_i\geq s+1.$
     The minimality of $s_i$ indicates that  \[w_p(\partial_L \mathcal{D}_{i}^{s})\leq w_p(\partial_L \mathcal{H}_{i}^{s})\leq \frac{\delta}{k},\]
     since otherwise $s$ would have been chosen as the smallest integer in step 7 of Algorithm~\ref{alg:hypergraph_container}.
     It remains to evaluate $w_p(\partial_L \mathcal{F}_{i}^{s})$. 
     If $L_i\subseteq I$, then \[w_p(\partial_L \mathcal{F}_{i}^{ s})=\mathbf{1}_{L\cap L_i=\emptyset}\cdot w_p(\partial_{L\cup L_i} \mathcal{H}_{i}^{s_i})\leq \frac{\delta}{k},\]
     where the last inequality is implied by the maximality of $L_i.$
     Otherwise, $L_i\not\subseteq I$. Then $ \mathcal{F}_{i}^{s}=\{L_i\}$ by definition by steps 9-13 in the algorithm, and 
    \[w_p(\partial_L \mathcal{F}_{i}^{s})\leq p^{s-|L|}\leq  p\leq \frac{\delta}{4k^2}\leq \frac{\delta}{k}.\]
  Combining both of the inequalities, we have $w_p(\partial_{L}\mathcal{D}_{i+1}^{s})\leq \delta/k+\delta/k= 2\delta/k.$
  Thus, the induction hypothesis holds for $i+1$, which completes our proof.
\end{proof}

Lemma~\ref{upper} immediately yields an  upper bound for $\mathcal {H}_i$ itself. 
Since for each $s$, the definition of $\mathcal{D}_i$ implies $\mathcal H_i^s \subseteq \mathcal H_0^s \cup \mathcal D_i^s$, taking
$\partial_L$ and applying Lemma~\ref{upper} gives the following corollary. 
\begin{cor}\label{upperH} For each $i$, and each subset $L\notin \langle \mathcal{H}_i\rangle$ and $s$ satisfying $1\leq |L|<s\leq k-1$, $$w_p(\partial_{L}\mathcal{H}_i^{s})\leq w_p(\partial_{L}\mathcal{H}_0^{s})+\frac{2\delta }{k}.$$
\end{cor}

\begin{lem}\label{smalls} Each output fingerprint
$|S|\leq 4k^2pn/\delta$.
\end{lem}
\begin{proof}
    We first analyze the edges added to or deleted from $\mathcal{D}_{i+1}$ at iteration $i$. 
    By Lemma~\ref{3.1,3.2,3.3}, every edge of $\mathcal{F}_i$ is newly added at iteration $i$, that is,
$\mathcal{F}_i \subseteq \mathcal{D}_{i+1} \setminus \mathcal{D}_i$.
On the other hand, when passing from $\mathcal{D}_i$ to $\mathcal{D}_{i+1}$, the only deleted
edges are those in $\langle \mathcal{F}_i \rangle \cap \mathcal{D}_i$. Therefore,
    \[w_p(\mathcal{D}_{i+1})-w_p(\mathcal{D}_i)\geq \sum_{E\in \mathcal{F}_i}p^{|E|}(1-w_p(\partial_E\mathcal{D}_{i})).\]
 Since $\mathcal{H}_0$ is $k$-bounded and the algorithm would not generate edges of size at least $k$, it naturally shows that $\mathcal{D}_i^{k}=\emptyset$ for each $i$. 
 Combining with Lemma~\ref{upper}, it implies that $w_p(\partial_{E}\mathcal{D}_i)\leq \sum_{s=2}^{k-1}w_p(\partial_{E}\mathcal{D}_i^s)+0\leq 2\delta.$
 Hence, 
    \[\sum_{E\in \mathcal{F}_i}p^{|E|}(1-w_p(\partial_E\mathcal{D}_{i}))\geq\sum_{E\in \mathcal{F}_i}p^{|E|}(1-2\delta)=(1-2\delta)w_p(\mathcal{F}_i)\geq \frac{w_p(\mathcal{F}_i)}{2},\]
    where the last inequality holds because $\delta<1/4$. 
    
Now consider those indices $i < J$ for which $L_i \subseteq I$.
In such a step, we add $L_i$ to the fingerprint, so
\[S_{i+1} = S_i \cup L_i,\qquad|L_i| \le k.\]
Thus each such step contributes at most $k$ new vertices to $S$.
Therefore the number of indices with $L_i \subseteq I$
is at least $|S|/k$. 
Moreover, whenever $L_i \subseteq I$, we have
$\mathcal{F}_i = \partial_{L_i} \mathcal{H}_i^{s_i}$ and the choice of $L_i$ implies $w_p(\partial_{L_i} \mathcal{H}_i^{s_i}) > \delta/k$.
As a consequence, we obtain the lower bound as
\begin{equation}\label{DJlow}
\begin{aligned}
w_p(\mathcal{D}_J) - w_p(\mathcal{D}_0)
&=
\sum_{i=0}^{J-1}
\bigl(w_p(\mathcal{D}_{i+1}) - w_p(\mathcal{D}_i)\bigr)
\ge
\sum_{\substack{0 \le i < J \\ L_i \subseteq I}}
\bigl(w_p(\mathcal{D}_{i+1}) - w_p(\mathcal{D}_i)\bigr)
\\
&\ge
\frac12
\sum_{\substack{0 \le i < J \\ L_i \subseteq I}}
w_p(\mathcal{F}_i)
>
\frac12
\cdot
\frac{|S|}{k}
\cdot
\frac{\delta}{k}
=
\frac{\delta |S|}{2k^2}.
\end{aligned}
\end{equation}

On the other hand, $\mathcal{D}_0 = \emptyset$, so $w_p(\mathcal{D}_0) = 0$.
By the stopping condition,
\[
w_p(\mathcal{H}_J^{>1}) \le \delta p |C_J|.
\]
All $1$-edges of $\mathcal{H}_J$ correspond to vertices outside $C_J$, so their
total weight is $p(2n - |C_J|)$.
Hence
\begin{equation}\label{DJup}
\begin{aligned}
w_p(\mathcal{D}_J)
\le
w_p(\mathcal{H}_J)
=
w_p(\mathcal{H}_J^{>1}) + p(2n - |C_J|)
\le
\delta p |C_J| + p(2n - |C_J|)
\le
2pn.
\end{aligned}
\end{equation}

Combining Eq.~\eqref{DJlow} and~\eqref{DJup} of $w_p(\mathcal{D}_J)$, we have $|S| \le 4k^2pn/\delta$.
\end{proof}

The next lemma estimates the time complexity of the algorithm itself.
It ensures that, for each input, the algorithm will cost polynomial time to get the output. 

\begin{lem}\label{timecomplex}
   The time complexity of Algorithm~\ref{alg:hypergraph_container}  is polynomial in $n$. 
\end{lem}
\begin{proof}
  Given a $k$-bounded hypergraph $\mathcal{H}$, an independent set $I$, constants $\lambda > 1,0<\delta <1/4$, we analyze the running time of the algorithm step by step.
  We begin with establishing an upper bound on $J$. 
  In each iteration $i$, the algorithm produces a non-empty edge set $\mathcal{F}_i$.
 For each $j>i$,  $\mathcal{F}_i\subseteq \langle \mathcal{H}_{j}\rangle$ since $\mathcal{F}_i\subseteq \mathcal H_{i+1}$ and $\langle\mathcal{H}_{i+1}\rangle\subseteq \langle\mathcal{H}_{j}\rangle$. 
  While, $\mathcal{F}_i\cap\langle\mathcal{H}_{i}\rangle=\emptyset$, which implies that $\mathcal{F}_i\cap \mathcal{F}_s=\emptyset$ for each $s<i.$
  Hence, $\mathcal{F}_i \cap(\cup_{s<i}\mathcal F_s)=\emptyset.$
  Moreover, the algorithm only generates edges of size less than $k$.
  Hence, $\cup_{s<J}\mathcal F_s \subseteq {[2n]\choose \leq (k-1)}.$
  It implies that 
  \[J\leq \sum_{s=0}^{J-1}|\mathcal{F}_s|= \left|\bigcup_{s=0}^{J-1}\mathcal{F}_s\right|\leq \sum_{i=1}^{k-1} \binom{2n}{i} = O(n^{k-1}).\] 
Once the iteration terminates at $J$, the next step verifies whether both $x_i$ and $\neg x_i$ are absent from $C_J$ for each $x_i$, which costs $O(n)$ time. 

  It suffices to examine the cost of each step within iteration $i$. 
  The algorithm computes $C_i$ by checking whether $\{v\}$ is an edge in $\mathcal{H}_i$ for each vertex $v$, which takes $O(n)$ time. 
  Computing $w_p(\mathcal{H}_i^{>1})$ requires summing $p^{|E|}$ over all edges $E$ of $\mathcal{H}_i^{>1}$. Since the number of edges is at most $O(n^k)$, this will cost $O(n^{k})$ time. 
  The search for a suitable pair $(s_i, L_i)$ requires $O(n^{2k})$ time by enumerating all $2\leq s\leq k$ and $|L|\leq k$, and scanning the edge set of $\mathcal{H}_i^{s}$. 
  The remaining operations (updating $\mathcal{F}_i$, $S_{i+1}$, and $\mathcal{H}_{i+1}$) involve only local modifications to the hypergraph, and thus the running time is polynomial in $n$. 
  Therefore, the total time complexity inside each iteration $i$ is polynomial in $n$. 
  Combining these bounds, we conclude that the algorithm runs in polynomial time in $n$.
\end{proof}

Specifically, we only need to input all independent sets with size less than $4k^2pn/\delta$ to generate and select all $C\in \mathcal{C}$. 
Checking whether a set with size less than $4k^2pn/\delta$ is independent will cost $O^*(2^{2H(2k^2p/\delta)n})$ time. 
This ensures that we only cost sub-exponential time to partition the original formula $\varphi$ into formulas $\varphi[C]$.
It remains to evaluate the time of solving $\varphi[C]$ for each $C\in \mathcal{C}.$
We next show that each container is itself not too large. 
The following lemma ensures that subsequent searches for satisfying assignments within a container can be performed efficiently, thereby controlling the overall running time. 

\begin{lem}\label{qsize}
   If $\mathcal{H}$ is a $(\lambda,p)_k$-structure with $p\leq 1/\lambda$, then each container $C\in \mathcal{C}$ satisfies $|C|
   \leq (2-1/\lambda)n$. 
\end{lem}

\begin{proof}
Recall that the maximal $p$-codegree $\Delta_{s,p}(\mathcal{H})$ is defined as $\Delta_{s,p}(\mathcal{H}) = \max_{|L| = s} \sum_{E \supseteq L} p^{|E|}$. 
This implies that $\Delta_{s,p}(\mathcal{H})\le \Delta_{2,p}(\mathcal{H})$  for all $2\leq s\leq k$.
Hence, we can estimate $\Delta_{s,p}$ as 
\begin{equation}
\Delta_{s,p}(\mathcal{H})\leq \Delta_{2,p}(\mathcal{H})\leq 
\lambda\bar{w}_p(\mathcal{H})p^{k}\leq \lambda\bar{w}_p(\mathcal{H})p^{s}\leq \bar{w}_p(\mathcal{H})p^{s-1}.  
\end{equation} The last inequality is due to $p\leq 1/\lambda$.

Note that in Algorithm~\ref{alg:hypergraph_container}, $\mathcal{H}_0=\mathcal{H}$.
We continue to evaluate the weight of edges inside $\mathcal{H}[C]$ via the algorithm. 
On the one hand, recall that $\langle\mathcal{H}_0\rangle\subsetneq \langle\mathcal{H}_{J}\rangle$ and each vertex in $C$ does not form a 1-edge. 
If an edge $F\in \mathcal{E}(\mathcal{H})$ satisfies $V(F)\subseteq C$, then there must exist an edge $E\in \mathcal{H}_J^{>1}$ such that $E\subseteq F$. 
This property shows that 
\begin{equation}\label{H0Cup}
\begin{aligned}
 w_p(\mathcal{H}[C])&\leq \sum_{s=2}^{k}\sum_{E\in \mathcal{G},|E|=s}\Delta_{s,p}(\mathcal{H})\leq \bar{w}_p(\mathcal{H})\left(\sum_{E\in \mathcal{G}}p^{|E|-1}\right)=\frac{w_p(\mathcal{H})}{2n}\cdot \frac{w_p(\mathcal{G})}{p}\\
 &\leq \delta \cdot w_p(\mathcal{H})< \frac{w_p(\mathcal{H})}{2},   
\end{aligned}
\end{equation}
where the last two inequalities hold because $w_p(\mathcal{G})\leq \delta p|C_J|\leq 2\delta pn$ and $\delta \leq 1/4$.
On the other hand, the lower bound on $w_p(\mathcal{H}[C])$ follows directly from the upper bound $\Delta_{1,p}(\mathcal{H})\leq \lambda\bar{w}_p(\mathcal{H})=\lambda w_p(\mathcal{H})/(2n)$, and further,
\begin{equation}\label{H0Clow}
\begin{aligned}
  w_p(\mathcal{H}[C])\geq w_p(\mathcal{H})-\Delta_{1,p}(\mathcal{H})(2n-|C|)\geq w_p(\mathcal{H})-\frac{\lambda\cdot w_p(\mathcal{H})}{2n}(2n-|C|).  
\end{aligned}
\end{equation}
Therefore, Eq.~\eqref{H0Cup} and~\eqref{H0Clow} imply $2n-|C|\geq n/\lambda$, which finishes the proof.
\end{proof}
We are ready to give a proof of Theorem~\ref{technical}.

\begin{proof}[Proof of Theorem~\ref{technical}]
We fix $p'=\min\{1/\lambda,\delta/(4k^2)\}$. 
Suppose that $\mathcal{H}_{\varphi}$ is the $(\lambda,p)_k$-structure with $p<p'$.
We construct families $\mathcal{S}$ and $\mathcal{C}$ in the following steps:
\begin{itemize}
    \item[(I)] Enumerate all independent sets of size less than $4k^2pn/\delta$ in $\mathcal{H}_{\varphi}$.
    \item[(II)] For each independent set $I$ with $|I|\leq 4k^2pn/\delta$, input $\mathcal{H}=\mathcal{H}_{\varphi},I,\delta$ and $p$.
    Run Algorithm~\ref{alg:hypergraph_container}  to receive $S,C$ and $\mathcal{G}$.
\item[(III)] Let $\mathcal{S}=\{S\mid S\neq $ \textsc{Error}$\}$, $\mathcal{C}=\{C\mid C\neq $ \textsc{Error}$\}$.
\end{itemize}
We show that these steps (I) and (II) are enough to generate all containers. 
By Lemma~\ref{3.11}, for different independent sets $I,I'$, $S=S'$ implies $C=C'$. 
Moreover, the rules for choosing $(s_i, L_i)$ in the algorithm imply that if $S$ is the output fingerprint set of $I$, then inputting $S$ itself will generate the same fingerprint. 
Hence, inputting all ``possible fingerprint sets'' is enough to find all possible containers. 
While Lemma~\ref{smalls} gives the upper bound on the size of fingerprints, it is enough to input all independent sets with size less than $4k^2pn/\delta$.

We now verify that sets $\mathcal{S}$ and $\mathcal{C}$ satisfy (1)-(4) in  Theorem~\ref{technical}.
By Lemma~\ref{3.11}, $I_{\alpha}\subseteq C_J$. 
Since $I_{\alpha}$ contains one literal for each variable, it implies that $C_J$ is not changed to \textsc{Error} in the algorithm. 
Hence, $C=C_J$ is added into $\mathcal{C}$. 
As a consequence, for each satisfying assignment $\alpha$ of $\varphi$, there exists $C\in \mathcal{C}$ such that $I_{\alpha}\subseteq C$, completing the proof of Theorem~\ref{technical} (1). 
By Lemma~\ref{smalls} every fingerprint produced by the algorithm has size \(\le 4k^2pn/\delta\), which implies Theorem~\ref{technical} (2). 
For Theorem~\ref{technical} (4),  we have \(|C|\le (2-1/\lambda)n\) by Lemma~\ref{qsize}. 
Moreover, the condition to stop the iteration implies \(w_p(\mathcal{G})\le \delta p|C|\le 2\delta pn\). 
By the same estimation as Eq.~\eqref{H0Cup}, it yields \(w_p(\mathcal{H}[C])\le \delta \cdot w_p(\mathcal{H})\).
Theorem~\ref{technical} (3) is naturally implied by steps 16-17 in Algorithm~\ref{alg:hypergraph_container}.  
It remains to show the time complexity. 
Step (I) requires us to enumerate all independent sets of size at most $4k^2pn/\delta.$ 
This can be done in time at most \({2n\choose \leq 4k^2pn/\delta}\leq 2^{2 H(2k^2p/\delta)n}\) . 
For each independent set $I$,  Lemma~\ref{timecomplex} implies that each run of the algorithm takes polynomial time in $n$. 
Hence, step (II) takes polynomial time.
Step (III) only requires us to count those $S_J, C_J$ that are not changed to \textsc{Error} and thus takes negligible additional time. 
So the total time is \(O^*(2^{2H(2k^2p/\delta)n})\). 
\end{proof}

Finally, we prove 
Corollary~\ref{regular} and~\ref{densemax}, where we employ a slightly different approach to bound the size of each container $C$. 
Recall that in the container algorithm, if $\mathcal{F}_i$ is 1-uniform, then $\mathcal{F}_i$ will bring new 1-edges to $\mathcal{H}_{i+1}$ and thus the size of $C_i$ decreases. 
Otherwise, the size of $C_i$ remains unchanged. 
The following lemma shows that when all elements of $\mathcal{F}_i$ have size greater than one, the weight $w_p(\mathcal{H}_i^{>1})$ does not decrease.
\begin{lem}\label{notde}
    Let $p=1/(2n)$ and $\mathcal{H}$ be a $k$-bounded hypergraph on $2n$ vertices. If $\mathcal{F}_i$ is $s$-uniform with $s\geq 2$, then
    \[w_p(\mathcal{H}_{i+1}^{>1})-w_p(\mathcal{H}_i^{>1})\geq 0.\]  
\end{lem}

\begin{proof}
 Similarly to the proof of Lemma~\ref{smalls}, we have 
\[
w_p(\mathcal{H}_{i+1}^{>1}) - w_p(\mathcal{H}_i^{>1}) \geq \sum_{E \in \mathcal{F}_i} \frac{1}{(2n)^s} \bigl(1 - w_p(\partial_E \mathcal{H}_i)\bigr).
\]
It remains to bound \(w_p(\partial_E \mathcal{H}_i)\). 
For set $E$ with $|E|=s<k$ and integer $0<j\leq k-s$, let \(a_{E, s+j}\) be the number of edges \(F \in \mathcal{H}_i\) such that \(E \subset F\) and \(|F| = s+j\); equivalently, this is the number of edges of size \(j\) in \(\partial_E \mathcal{H}_i\).
Since \(\mathcal{H}_i\) is an antichain by Lemma~\ref{3.1,3.2,3.3}, it follows that \(\partial_E \mathcal{H}_i\) is also an antichain. 
Hence, the LYM inequality (Lemma~\ref{LYM}) bounds the number of edges inside \(\partial_E \mathcal{H}_i\) as 
\[
\sum_{j=1}^{k-s} \frac{a_{E,s+j}}{\binom{2n-|E|}{j}} \le 1.
\]
Observing that $\binom{2n-s}{j} \le (2n)^j$ for every $j$, we have $$\frac{a_{E,s+j}}{(2n)^j} \le \frac{a_{E,s+j}}{\binom{2n-s}{j}}.$$ This yields the upper bound of $w_p(\partial_E \mathcal{H}_i)$, which is 
\[
w_p(\partial_E \mathcal{H}_i) = \sum_{j=1}^{k-s} \frac{a_{E,s+j}}{(2n)^j} \le \sum_{j=1}^{k-s} \frac{a_{E,s+j}}{\binom{2n-s}{j}} \le 1.
\]
Therefore, $1-w_p(\partial_E \mathcal{H}_i)\geq 0$ for each $E\in\mathcal{F}_i$.
This is enough to obtain 
$w_p(\mathcal{H}_{i+1}^{>1}) - w_p(\mathcal{H}_i^{>1}) \ge 0$.
\end{proof}

\begin{proof}[Proof of Corollary~\ref{regular}]
   We select the parameters as $p=1/(2n), \delta=d/3$ and run Algorithm~\ref{alg:hypergraph_container} for $\mathcal{H}$.
   We first verify that these parameters satisfy the assumptions of the algorithm.
Since every edge of $\mathcal H$ has size at least $2$ and $\mathcal{H}$ is an antichain, its edges satisfy 
\[d\leq \sum_{s=2}^k\frac{|\mathcal{H}^s|}{(2n)^s}\leq \frac{1}{2}\left(\sum_{s=2}^{k}|\mathcal{H}^s|{2n \choose s}^{-1}\right)\leq \frac{1}{2},\]
where the second inequality holds because $s!\geq 2$, and the last inequality is implied by Lemma~\ref{LYM}. Therefore,
$\delta=d/3 \le 1/6 < 1/4$.
Moreover, $p=1/(2n)\le \delta/(4k^2)$ for sufficiently large $n$.
Thus, all requirements of the algorithm are satisfied.
We obtain $\mathcal{C}$ by enumerating all possible fingerprints as in Theorem~\ref{technical}.
This takes time $O^*({2n\choose \leq 6k^2/d})$, which is polynomial in $n$.

Subsequently, it remains to calculate the time complexity of searching for satisfying assignments inside each reduced formula $\varphi[C]$. 
  
\begin{claim}
  $|C|\leq 2n-d/(2\Delta_{1,p}(\mathcal{H}))$. 
\end{claim}
  \begin{proof}
    Let $\lambda=2n \Delta_{1,p}(\mathcal{H})/d$.
    For iteration $i$, suppose that $\mathcal{F}_i$ is $s$-uniform for some $1<s\leq k-1$, then Lemma~\ref{notde} yields 
$$w_p(\mathcal{H}_{i+1}^{>1})-w_p(\mathcal{H}_i^{>1})\geq 0, $$
and the size of $C_i$ does not change.
Otherwise, $\mathcal{F}_i$ is $1$-uniform. 
Let $q_i = |\mathcal{F}_i|$. 
Equivalently, $q_i$ denotes the number of vertices deleted from $C_i$ for each iteration $i$.
In this case, we have 
\[w_p(\mathcal{H}_{i+1}^{>1})-w_p(\mathcal{H}_i^{>1})\geq \sum_{\{v\}\in\mathcal{F}_i}-pw_p(\partial_v\mathcal{H}_{i})\geq -\sum_{\{v\}\in\mathcal{F}_i}p(w_p(\partial_v\mathcal{H}_0)+2\delta)\geq -(2\delta +d\lambda)pq_i.\]
The second inequality is implied by Lemma~\ref{upperH}, and the last inequality follows from $pw_p(\partial_v\mathcal{H}_0)\leq \Delta_{1,p}(\mathcal{H})=\lambda dp$. Consequently, 
\begin{equation}\label{H0>1low}
\begin{aligned}
w_p(\mathcal{H}_J^{>1})-w_p(\mathcal{H}_0^{>1})\geq \sum_{i=0}^{J-1}(w_p(\mathcal{H}_{i+1}^{>1})-w_p(\mathcal{H}_i^{>1}))\geq -(2\delta +d\lambda)pq,      
\end{aligned}
\end{equation}
where $q = \sum_i q_i$ is the total number of vertices out of $C$. 
On the other hand, the stopping condition of the algorithm requires $w_p(\mathcal{H}_J^{>1})\leq \delta p|C|$. 
The fact that $\mathcal{H}_0$ does not contain 1-edge implies $w_p(\mathcal{H}_0^{>1})=w_p(\mathcal{H}_0).$
Thus,
\begin{equation}\label{H0>1up}
\begin{aligned}
w_p(\mathcal{H}_J^{>1})-w_p(\mathcal{H}_0^{>1})= w_p(\mathcal{H}_J^{>1})-w_p(\mathcal{H}_0) \leq p\delta(2n-q)-d=\delta-d-\delta pq.   
\end{aligned}
\end{equation}
Combining Eq.~\eqref{H0>1low} and~\eqref{H0>1up} gives
\[(2\delta +d\lambda)pq\geq d-\delta+\delta pq, \]
which simplifies to 
\[q\geq \frac{d-\delta}{(\delta +d\lambda)p}.\]
By $p=1/(2n),\delta=d/3, \lambda=2n \Delta_{1,p}(\mathcal{H})/d$, we have 
\[\frac{d-\delta}{(\delta +d\lambda)p}=\frac{4n}{1+3\lambda}\geq \frac{n}{\lambda}=\frac{d}{2\Delta_{1,p}(\mathcal{H})},\]
and thus each container $C$ satisfies $|C|\leq 2n-d/(2\Delta_{1,p}(\mathcal{H}))$.  
  \end{proof}
By the claim, applying the assumed $O^*(\beta^n)$ algorithm to each reduced formula
$\varphi[C]$ takes $O^*(\beta^{n-d/(2\Delta_{1,p}(\mathcal{H}))})$ time. 
Since the number of containers is polynomial,
the time to solve all $\varphi[C]$s is $O^*(\beta^{n-d/(2\Delta_{1,p}(\mathcal{H}))})$. 
This implies that the total running time is $O^*(\beta^{n-d/(2\Delta_{1,p}(\mathcal{H}))})$. 
\end{proof}

\begin{proof}[Proof of Corollary~\ref{densemax}]
For sufficiently large $n$, we input $\mathcal{H}=\mathcal{H}_{\varphi},p=1/(2n),0<\delta'< \min\{1/4,\delta d\}$ and run Algorithm~\ref{alg:hypergraph_container} to obtain $\mathcal{C}$. 
If $\mathcal{C}=\emptyset$, then it implies that the SAT formula is unsatisfiable. 
Otherwise, $\mathcal{C}\neq \emptyset$.
For each $E \in \mathcal{G}$, if $E\notin \mathcal{H}$, then the $p$-weight of edges in $\mathcal{H}$ containing $E$ can be bounded as
\[
\sum_{\substack{F \in \mathcal{H} \\ E \subseteq F}} \frac{1}{(2n)^{|F|}}
= \sum_{s =|E|+1}^{k} \frac{a_{E,s}}{(2n)^s}
= \frac{1}{(2n)^{|E|}} \left( \sum_{s = |E|+1}^{k} \frac{a_{E,s}}{(2n)^{s-|E|}} \right)
\leq  \frac{1}{(2n)^{|E|}} \left( \sum_{s = |E|+1}^{k} \frac{a_{E,s}}{\binom{2n}{s-|E|}} \right)
\leq \frac{1}{(2n)^{|E|}},
\]
where the last inequality follows from Lemma~\ref{LYM}.
Otherwise, $E\in \mathcal{G}\cap \mathcal{H}$. Since $\mathcal{H}$ is an antichain, the $p$-weight of edges in $\mathcal{H}$ containing $E$ is $1/(2n)^{|E|}.$
Consequently, after Algorithm~\ref{alg:hypergraph_container} terminates, for each container $C\in \mathcal{C}$, the hypergraph $\mathcal{H}[C]$ has $p$-weight as
\[
w_p(\mathcal{H}[C])
\leq \sum_{E \in \mathcal{G}} \frac{1}{(2n)^{|E|}}
= w_p(\mathcal{G})
\leq 2\delta' p n< 2\delta d pn=\delta d.
\]
Implied by the condition $w_p(\mathcal{H})\geq d$, we have $w_p(\mathcal{H}[C])< \delta d\leq \delta \cdot w_p(\mathcal{H})$ . 
This shows that for each assignment $\alpha$, if its corresponding set is included by some container $C$, then the weight of the unsatisfied clauses is at most $\delta d\leq\delta w_p(\mathcal{H}_{\varphi})$.
The time complexity is determined by the container algorithm, which is $O^*((2n)^{2k^2/\delta'}) = O^*(1)$.

\end{proof}

\section{Conclusion and Future Directions}\label{sec5}

In this work, we extend the hypergraph container method to general SAT formulas by introducing the notion of  $(\lambda,p)_k$-structures, as a generalization of $(D, c, \epsilon)$-structures, which allows us to derive algorithmic consequences for approximate MAX-E$k$-SAT and general SAT.
By weighting clauses according to their lengths and adapting the container construction to this weighted non-uniform setting, we obtain an explicit deterministic procedure that covers all satisfying assignments by a small family of containers. We conclude by highlighting two directions for future work.

\paragraph{Further Algorithmic Applications of Theorem~\ref{technical}}
Theorem~\ref{technical} serves as a tailored version of the container method in SAT, offering a new perspective for both approximation algorithms and exact SAT solving. 
Beyond these immediate consequences, the framework suggests several broader applications.

\textbf{(1) CSP Property Testing: } The $\varepsilon$-far SAT formulas require the removal of at least $\varepsilon \binom{n}{2}$ clauses to become satisfiable. 
Such formulas naturally satisfy the globally dense structure, for which the speedup provided by the container method is the most efficient.
This raises a question for CSP property testing: can Theorem~\ref{technical} simplify the container-based proofs of Blais and Seth~\cite{blais2024new}, or improve their sample complexity bounds? 

\textbf{(2) SAT Solution Counting: } The \#SAT problem, which asks for the number of satisfying assignments of a SAT formula, has been extensively studied~\cite{valiant1979complexity,zhou2010new,kutzkov2007new,dahllof2005counting,peng2025new}. 
A natural question is whether Theorem~\ref{technical} can yield a similar acceleration for \#SAT, pushing the current algorithmic limits.

\paragraph{Applications in Practical SAT Solving}
While the container method originates in extremal combinatorics, its core intuition may also benefit applied SAT solving. 
Modern SAT solvers, particularly those relying on Conflict-Driven Clause Learning (CDCL)~\cite{biere2009handbook}, often face exponential bottlenecks when dealing with dense, globally entangled formulas. 
The container algorithm acts as a powerful space-partitioning tool, aligning with the Cube-and-Conquer paradigm~\cite{heule2011cube}. 
Just as Cube-and-Conquer partitions hard formulas into independent subproblems by generating cubes, the container method uses fingerprint enumeration as a preprocessing step to reduce the global problem into subproblems of much smaller size.
This reduces the global complexity to manageable local search spaces.

\section*{Acknowledgments}
This research was supported by the National Key Research and Development Program of China (Grant No. 2023YFA1010201 and 2020YFA0713100), the National Natural Science Foundation of China (Grant No. 12125106, 12171452, 12231014 and 125B2019), and the Innovation Program for Quantum Science and Technology (Grant No. 2021ZD0302902). 

We thank Cameron Seth for helpful discussions on the container method and its applications. We thank Marcelo Campos and Wojciech Samotij for their helpful comments, and Yihan Chen for helpful discussions and valuable suggestions.

\bibliographystyle{alpha} 
\bibliography{refs}

@inproceedings{balogh2023nearly,
  title     = {Nearly all $k$-{SAT} functions are unate},
  author    = {Balogh, J{\'o}zsef and Dong, Dingding and Lidick{\'y}, Bernard and Mani, Nitya and Zhao, Yufei},
  booktitle = {ACM Symposium on Theory of Computing (STOC)},
  year      = {2023},
  organization = {ACM}
}

@article{balogh2015independent,
  title={Independent sets in hypergraphs},
  author={Balogh, J{\'o}zsef and Morris, Robert and Samotij, Wojciech},
  journal={Journal of the American Mathematical Society},
  volume={28},
  number={3},
  pages={669--709},
  year={2015}
}

@article{balogh2020efficient,
  title={An Efficient Container Lemma},
  author={Balogh, J{\'o}zsef and Samotij, Wojciech},
  journal={Discrete Analysis},
  volume={2020},
  pages={Paper No. 17, 56 pp.},
  year={2020},
  doi={10.19086/da.17354}
}

@book{biere2009handbook,
  title={Handbook of Satisfiability},
  editor={Biere, Armin and Heule, Marijn and van Maaren, Hans and Walsh, Toby},
  volume={185},
  year={2009},
  publisher={IOS Press}
}

@inproceedings{blais2024new,
  title={New graph and hypergraph container lemmas with applications in property testing},
  author={Blais, Eric and Seth, Cameron},
  booktitle={Proceedings of the 56th Annual ACM Symposium on Theory of Computing},
  pages={1793--1804},
  year={2024}
}

@article{blais2025testing,
  title={Testing graph properties with the container method},
  author={Blais, Eric and Seth, Cameron},
  journal={SIAM Journal on Computing},
  volume={54},
  number={6},
  pages={1626--1642},
  year={2025},
  publisher={SIAM}
}

@article{bollobas1965generalized,
  title     = {On generalized graphs},
  author    = {Bollob{\'{a}}s, B.},
  journal   = {Acta Mathematica Academiae Scientiarum Hungaricae},
  volume    = {16},
  pages     = {447--452},
  year      = {1965}
}

@article{campos2024towards,
  title={Towards an optimal hypergraph container lemma},
  author={Campos, Marcelo and Samotij, Wojciech},
  journal={arXiv preprint arXiv:2408.06617},
  year={2024}
}

@inproceedings{dong2022enumerating,
  title     = {Enumerating $k$-{SAT} functions},
  author    = {Dong, Dingding and Mani, Nitya and Zhao, Yufei},
  booktitle = {ACM-SIAM Symposium on Discrete Algorithms (SODA)},
  year      = {2022},
  pages = {2141--2184},
  organization = {ACM and SIAM}
}

@inproceedings{hansen2019faster,
  title={Faster k-sat algorithms using biased-ppsz},
  author={Hansen, Thomas Dueholm and Kaplan, Haim and Zamir, Or and Zwick, Uri},
  booktitle={Proceedings of the 51st Annual ACM SIGACT Symposium on Theory of Computing},
  pages={578--589},
  year={2019}
}

@article{haastad2001some,
  title={Some optimal inapproximability results},
  author={H{\aa}stad, Johan},
  journal={Journal of the ACM (JACM)},
  volume={48},
  number={4},
  pages={798--859},
  year={2001},
  publisher={ACM New York, NY, USA}
}

@article{hertli2014,
  title={3-SAT faster and simpler---Unique-SAT bounds for PPSZ hold in general},
  author={Hertli, Timon},
  journal={SIAM Journal on Computing},
  volume={43},
  number={2},
  pages={718--729},
  year={2014},
  publisher={SIAM}
}

@inproceedings{heule2011cube,
  title={Cube and conquer: Guiding CDCL SAT solvers by lookaheads},
  author={Heule, Marijn JH and Kullmann, Oliver and Wieringa, Siert and Biere, Armin},
  booktitle={Hardware and Software: Verification and Testing: 7th International Haifa Verification Conference, HVC 2011, Haifa, Israel, December 6-8, 2011, Revised Selected Papers 7},
  pages={50--65},
  year={2012},
  organization={Springer}
}

@article{impagliazzo2001complexity,
  title={On the complexity of k-SAT},
  author={Impagliazzo, Russell and Paturi, Ramamohan},
  journal={Journal of Computer and System Sciences},
  volume={62},
  number={2},
  pages={367--375},
  year={2001},
  publisher={Elsevier}
}

@article{impagliazzo2001problems,
  title={Which problems have strongly exponential complexity?},
  author={Impagliazzo, Russell and Paturi, Ramamohan and Zane, Francis},
  journal={Journal of computer and system sciences},
  volume={63},
  number={4},
  pages={512--530},
  year={2001},
  publisher={Elsevier}
}

@article{kleitman1982number,
  title={On the number of graphs without 4-cycles},
  author={Kleitman, Daniel J and Winston, Kenneth J},
  journal={Discrete Mathematics},
  volume={41},
  number={2},
  pages={167--172},
  year={1982},
  publisher={Elsevier}
}

@article{lubell1966short,
  title     = {A short proof of Sperner's lemma},
  author    = {Lubell, D.},
  journal   = {Journal of Combinatorial Theory},
  volume    = {1},
  pages     = {299},
  year      = {1966}
}

@article{mesalkin1963generalization,
  title     = {A generalization of {Sperner}'s theorem on the number of subsets of a finite set},
  author    = {Me{\v{s}}alkin, L. D.},
  journal   = {Teoriya Veroyatnostei i ee Primeneniya},
  volume    = {8},
  pages     = {219--220},
  year      = {1963}
}

@article{nenadov2024probabilistic,
  title={Probabilistic hypergraph containers},
  author={Nenadov, Rajko},
  journal={Israel Journal of Mathematics},
  volume={261},
  number={2},
  pages={879--897},
  year={2024},
  publisher={Springer}
}

@article{paturi2005improved,
  title={An improved exponential-time algorithm for k-SAT},
  author={Paturi, Ramamohan and Pudl{\'a}k, Pavel and Saks, Michael E and Zane, Francis},
  journal={Journal of the ACM (JACM)},
  volume={52},
  number={3},
  pages={337--364},
  year={2005},
  publisher={ACM New York, NY, USA}
}

@article{saxton2015hypergraph,
  title={Hypergraph containers},
  author={Saxton, David and Thomason, Andrew},
  journal={Inventiones mathematicae},
  volume={201},
  number={3},
  pages={925--992},
  year={2015},
  publisher={Springer}
}

@article{saxton2016simple,
  title={Simple containers for simple hypergraphs},
  author={Saxton, David and Thomason, Andrew},
  journal={Combinatorics, Probability and Computing},
  volume={25},
  number={3},
  pages={448--459},
  year={2016},
  publisher={Cambridge University Press}
}

@inproceedings{scheder2021ppsz,
  title={PPSZ is better than you think},
  author={Scheder, Dominik},
  booktitle={2021 IEEE 62nd Annual Symposium on Foundations of Computer Science (FOCS)},
  pages={205--216},
  year={2021},
  organization={IEEE},
  doi={10.1109/FOCS52979.2021.00028}
}

@inproceedings{seth2025tolerant,
  title={A Tolerant Independent Set Tester},
  author={Seth, Cameron},
  booktitle={Proceedings of the 57th Annual ACM Symposium on Theory of Computing},
  pages={1031--1042},
  year={2025}
}

@article{yamamoto1954logarithmic,
  title     = {Logarithmic order of free distributive lattice},
  author    = {Yamamoto, K.},
  journal   = {Journal of the Mathematical Society of Japan},
  volume    = {6},
  pages     = {343--353},
  year      = {1954}
}

@inproceedings{zamir2023algorithmic,
  title={Algorithmic Applications of Hypergraph and Partition Containers},
  author={Zamir, Or},
  booktitle={Proceedings of the 55th Annual ACM Symposium on Theory of Computing (STOC)},
  pages={985--998},
  year={2023},
  publisher={ACM}
}

@inproceedings{arora1995polynomial,
  title={Polynomial time approximation schemes for dense instances of NP-hard problems},
  author={Arora, Sanjeev and Karger, David and Karpinski, Marek},
  booktitle={Proceedings of the twenty-seventh annual ACM symposium on Theory of computing},
  pages={284--293},
  year={1995}
}

@article{valiant1979complexity,
  title={The complexity of computing the permanent},
  author={Valiant, Leslie G},
  journal={Theoretical Computer Science},
  volume={8},
  number={2},
  pages={189--201},
  year={1979},
  publisher={Elsevier},
  doi={10.1016/0304-3975(79)90043-0}
}

@inproceedings{zhou2010new,
  title={New worst-case upper bound for \#2-SAT and \#3-SAT with the number of clauses as the parameter},
  author={Zhou, Junping and Yin, Minghao and Zhou, Chunguang},
  booktitle={Proceedings of the Twenty-Fourth AAAI Conference on Artificial Intelligence (AAAI)},
  pages={217--222},
  year={2010},
  editor={Fox, Maria and Poole, David},
  publisher={AAAI Press},
  doi={10.1609/aaai.v24i1.7537}
}

@article{kutzkov2007new,
  title={New upper bound for the \#3-SAT problem},
  author={Kutzkov, Konstantin},
  journal={Information Processing Letters},
  volume={105},
  number={1},
  pages={1--5},
  year={2007},
  publisher={Elsevier},
  doi={10.1016/j.ipl.2007.07.014}
}

@article{dahllof2005counting,
  title={Counting models for 2SAT and 3SAT formulae},
  author={Dahll{\"o}f, Vilhelm and Jonsson, Peter and Wahlstr{\"o}m, Magnus},
  journal={Theoretical Computer Science},
  volume={332},
  number={1-3},
  pages={265--291},
  year={2005},
  publisher={Elsevier},
  doi={10.1016/j.tcs.2004.10.022}
}

@inproceedings{peng2025new,
  title={New Algorithms for \#2-SAT and \#3-SAT},
  author={Peng, Junqiang and Sheng, Zimo and Xiao, Mingyu},
  booktitle={Proceedings of the Thirty-Fourth International Joint Conference on Artificial Intelligence (IJCAI)},
  year={2025},
  note={arXiv preprint arXiv:2507.14504}
}

\textit{\large Emails:}
\begin{itemize}
    \item \textit{hzcqj2020@mail.ustc.edu.cn}
    \item \textit{lyp\_inustc@mail.ustc.edu.cn}
    \item \textit{jiema@ustc.edu.cn}
    \item \textit{drzhangx@ustc.edu.cn}
\end{itemize}
\newpage
\appendix
\section{Appendix}\label{Appendix}

\begin{recallprop}{bigrandom}
For each \(k\ge3\), \(p\in(0,1]\), and \(\lambda>k\):
\begin{itemize}
    \item[(1)] There exists a constant
\(C_0=C_0(k,p,\lambda)>0\) such that if \( d\ge C_0\log n,\)
then \(\mathcal H_\Phi\) is a \((\lambda,p)_k\)-structure w.h.p.
as \(n\to\infty\).
\item[(2)] If  $d= o(\log n)$, then $\mathcal{H}_{\Phi}$ is not a
\((\lambda,p)_k\)-structure w.h.p. as \(n\to\infty\). 
\end{itemize}   \end{recallprop}

\begin{proof}
 For (1), let $C_0=\max \{1/(kh(\lambda/k))+1, 4/(\lambda p^k) +1\}$ and such $C_0$ will satisfy all conditions. 
 We verify the
three conditions of $(\lambda,p)_k$. 
Since \(\mathcal H_{\Phi}\) is obtained by
choosing \(k\)-edges uniformly from all
\(k\)-sets of literals containing no complementary pair, for each vertex $v$, its degree $\deg_{\mathcal H_{\Phi}}(v)\sim\text{Hyp}(2^k{n\choose k}, 2^{k-1}{n-1\choose k-1},m).$ The expectation of $\deg_{\mathcal H_{\Phi}}(v)$ is
\[\mathbb E [\deg_{\mathcal H_{\Phi}}(v)]=m\cdot\frac{2^{k-1}\binom{n-1}{k-1}}{2^k\binom nk}=m\cdot \frac{k}{2n}=kd.\]
By Chernoff bound, 
\[\Pr[\deg_{\mathcal H_{\Phi}}(v) > \lambda d]
\le
e^{-kd(h(\frac{\lambda}{k}))},   
\]
where $h(x)=x\ln{x}-x+1$. 
Hence
$\Pr[\Delta_1(\mathcal H_{\Phi}) > \lambda d]\le 2n\cdot \Pr[\deg_{\mathcal H_{\Phi}}(v) > \lambda d]\leq 2ne^{-kdh(\frac{\lambda}{k})}$.
Since \(\lambda>k\), \(h(\frac{\lambda}{k}) > 0\), and $2ne^{-kdh(\frac{\lambda}{k})}\le2n\cdot (\frac{1}{n})^{kh(\frac{\lambda}{k})C_0}$.
If $C_0> 1/(kh(\frac{\lambda}{k})), $ then $2ne^{-kdh(\frac{\lambda}{k})}$ tends to \(0\) as $n\to \infty$. 
Therefore,
$\Delta_1(\mathcal H_{\Phi})\le \lambda d$ with high probability.\\
\indent To clarify the condition of $\Delta_2$, we define 
$\operatorname{codeg}_{\mathcal H_{\Phi}}(u,v)
=
|\{E\in \mathcal{E}(\mathcal H_{\Phi}): \{u,v\}\subseteq E\}|$ for distinct literals $u,v$ with $u\not= \neg v$.
The expectation of $\operatorname{codeg}_{\mathcal H_{\Phi}}(u,v)$ is 
\[
\mathbb E[\operatorname{codeg}_{\mathcal H_{\Phi}}(u,v)]
=
m\cdot
\frac{2^{k-2}\binom{n-2}{k-2}}
{2^k\binom nk}
=
m\cdot \frac{k(k-1)}{4n(n-1)}
=
\frac{k(k-1)}{2(n-1)}d.
\]
\indent
By Chernoff bound, 
\[\Pr[\operatorname{codeg}_{\mathcal H_{\Phi}}(u,v) > \lambda dp^k ] \le \exp\left(-\frac{k(k-1)}{2(n-1)}d h\left(\frac{2\lambda p^k (n-1)}{k(k-1)}\right)\right).\]
When $n$ is sufficiently large, $h(\frac{2\lambda p^k(n-1)}{k(k-1)})\geq \frac{1}{2}\cdot \frac{2\lambda p^k(n-1)}{k(k-1)}(\log n+\log(\lambda p^k)-\log(k(k-1)))\geq \frac{\lambda p^k(n-1)}{2k(k-1)}\log n.$
Hence, when $n\to \infty$, the probability to have a large $\Delta_2$ is 
$\Pr[\Delta_2(\mathcal H_{\Phi}) > \lambda dp^k] \le 4 \binom{n}{2}\exp\left(-\frac{\lambda p^kd\log n}{4}\right)$.
Since $d\geq C_0\log n\geq (4/(\lambda p^k)+1)\log n,$
\[
\Pr[\Delta_2(\mathcal H_\Phi)>\lambda d p^k]
\le
4\binom n2
\exp\left(
-\left(1+\frac{\lambda p^k}{4}\right)(\log n)^2
\right)
\to 0.
\]
\indent
Finally, the condition $d\ge p^{1-k}$ holds as $n\to \infty$ since \(p\) is fixed and \(d\ge C_0\log n\).
This completes the proof of (1).\\
\indent 
For (2), we assume that $d\geq p^{1-k}$ and $\lambda d\geq 1$, since otherwise the $(\lambda,p)_k$ condition directly fails. 
In this case, it suffices to prove that $\Delta_1(\mathcal{H}_{\Phi})> \lambda d$ happens with high probability as $n\to \infty.$
For $t=\lfloor \lambda d\rfloor+1,$
we define random variable \(X=\bigl|\{v\in V(\mathcal H_\Phi):\deg_{\mathcal H_\Phi}(v)\ge t\}\bigr|\) be the number of vertices with degree at least $t$.
For each $v$, \(\deg_{\mathcal H_\Phi}(v)\sim \operatorname{Hyp}(N,K,m)\), where \(N=2^k\binom{n}{k}\) and \(K=2^{k-1}\binom{n-1}{k-1}\). 
Hence
\[
\Pr[\deg_{\mathcal H_\Phi}(v)=t]
=
\frac{\binom Kt\binom{N-K}{m-t}}{\binom Nm}.
\]
Since \(d=o(\log n)\), we have \(t=O(d)=o(K)\), \(m=o(N)\), and \(mK/N=kd\). Therefore, 
\[
\Pr[\deg_{\mathcal H_\Phi}(v)=t]
\ge
\exp(-O(d))\binom mt\left(\frac KN\right)^t
=
\exp(-O(d))\binom{2dn}{t}\left(\frac{k}{2n}\right)^t.
\]
By Stirling's formula, \(t!=\Theta(\sqrt t (t/e)^t)\) implies \(\frac{(kd)^t}{t!}
=
\Theta\left(
\frac1{\sqrt t}
\left(\frac{ekd}{t}\right)^t
\right).
\)
Since \(t=\lambda d+O(1)\), the right-hand side is at least \(\exp(-O(d))\). Therefore
\[
\exp(-O(d))\binom{2dn}{t}\left(\frac{k}{2n}\right)^t
\ge
\exp(-O(d))\frac{(kd)^t}{t!}
\ge
\exp(-O(d)).
\]
Thus there exists \(C_1=C_1(k,\lambda)>0\) such that
$\Pr[\deg_{\mathcal H_\Phi}(v)=t]\ge \exp(-C_1d).$
Moreover, 
\[
\mathbb E[X]
\ge
2n\Pr[\deg_{\mathcal H_\Phi}(v)=t]
\ge
2n\exp(-C_1d)=2n^{(1-o(1))},
\]
where the last equality follows from \(\exp(-C_1d)=n^{-o(1)}\). This implies that \(\mathbb E[X]\to\infty\).

 To apply Chebyshev's inequality, we need to show $\text{Var}(X)=o(\mathbb{E}[X]^2)$. Let $X_u=\mathbf{1}_{\deg (u)\geq t}$. 
 Then \[\mathbb{E}[X^2]=\mathbb{E}\left[\left(\sum_{u}X_u\right)^2\right]=\sum_{u}\mathbb{E}[X_u]+\sum_{u\not=v}\mathbb{E}[X_uX_v]=\mathbb{E}[X]+\sum_{u\not=v}\mathbb{E}[X_uX_v],\]
 where the second equality holds because $X_u^2=X_u$ for each $u$. 
 Since $\mathbb E[X] \to \infty$, it remains to prove 
\[
\Pr\Bigl[
\deg_{\mathcal H_\Phi}(u)\ge t,\,
\deg_{\mathcal H_\Phi}(v)\ge t
\Bigr]
\le
(1+o(1))
\Pr\bigl[\deg_{\mathcal H_\Phi}(u)\ge t\bigr]
\Pr\bigl[\deg_{\mathcal H_\Phi}(v)\ge t\bigr],
\]
for all distinct literals \(u,v\). 
Let $J$ be the number of chosen edges containing both $u,v.$
If $u=\neg v$, then $\Pr[J\geq 1]=0=o\left(\Pr[\deg_{\mathcal H_\Phi}(u)\geq t]^2\right)$.
Otherwise, the number of edges containing both $u$ and $v$ is $2^{k-2}{n-2\choose k-2}$.
It implies $J\sim \text{Hyp}(2^k{n\choose k}, 2^{k-2}{n-2\choose k-2}, m)$ and $\mathbb{E}(J)=\frac{k(k-1)m}{4n(n-1)}=o(\log n/n)$.  
Hence, \[\Pr[J\geq 1]\leq \mathbb{E}(J)=o\left(\frac{\log n}{n}\right)=o(\exp(-2C_1d))=o\left(\Pr[\deg_{\mathcal H_\Phi}(u)=t]^2 \right) = o\left(\Pr[\deg_{\mathcal H_\Phi}(u)\geq t]^2\right).\]
As a consequence, $\Pr[J=0]=1-o\left(\Pr[\deg_{\mathcal H_\Phi}(u)\geq t]^2\right)$ in both cases.

It suffices to estimate the probability of the event ``both $\deg (u)\geq t$ and $\deg (v)\geq t$'' conditioned on the event \(A=\{J=0\}\). 
Assume that \(J=0\) and \(\deg_{\mathcal H_\Phi}(u)=x\) happen. After exposing the \(x\) chosen edges containing \(u\), the remaining \(m-x\) chosen edges are uniformly sampled from the admissible edges not containing \(u\).
Therefore,
\[
\deg_{\mathcal H_\Phi}(v)
\mid
\bigl(
\deg_{\mathcal H_\Phi}(u)=x,\ A
\bigr)
\sim
\left\{
\begin{aligned}
 &\operatorname{Hyp}\left(
N-K,\,
K-2^{k-2}\binom{n-2}{k-2},\,
m-x
\right),\qquad &u\not = \neg v,\\   &\operatorname{Hyp}\left(
N-K,\,
K,\,
m-x
\right),\qquad &u = \neg v.
\end{aligned}
\right.
\]
Since $m-x$ decreases as $x$ increases, the probability $\Pr[\deg_{\mathcal H_\Phi}(v)\geq t | \deg_{\mathcal H_\Phi}(u)=x,A]$ is nonincreasing in $x$.
Hence,
\begin{equation*}
\begin{aligned}
\Pr\!\left[
\deg_{\mathcal H_\Phi}(u)\ge t \,\middle|\, A
\right]
&=
\Pr\!\left[
\deg_{\mathcal H_\Phi}(v)\ge t \,\middle|\, A
\right]  \\
&=
\sum_{x=0}^{t-1}
\Pr\!\left[
\deg_{\mathcal H_\Phi}(v)\ge t
\,\middle|\,
\deg_{\mathcal H_\Phi}(u)=x,\ A
\right]
\Pr\!\left[
\deg_{\mathcal H_\Phi}(u)=x
\,\middle|\,
A
\right]                                      \\
&\quad+
\Pr\!\left[
\deg_{\mathcal H_\Phi}(v)\ge t
\,\middle|\,
\deg_{\mathcal H_\Phi}(u)\ge t,\ A
\right]
\Pr\!\left[
\deg_{\mathcal H_\Phi}(u)\ge t
\,\middle|\,
A
\right]                                      \\
&\ge
\Pr\!\left[
\deg_{\mathcal H_\Phi}(v)\ge t
\,\middle|\,
\deg_{\mathcal H_\Phi}(u)\ge t,\ A
\right].
\end{aligned}
\end{equation*} 
This implies that \begin{equation*}
\begin{aligned}
\Pr[\deg_{\mathcal H_\Phi}(u)\geq t,\deg_{\mathcal H_\Phi}(v)\geq t, A] &\leq \frac{\Pr[\deg_{\mathcal H_\Phi}(u)\geq t,A]^2}{\Pr[A]}\leq \frac{\Pr[\deg_{\mathcal H_\Phi}(u)\geq t]^2}{\Pr[A]}\\
&=(1+o(1))\Pr[\deg_{\mathcal H_\Phi}(u)\geq t]^2.    
\end{aligned}
\end{equation*}
By combining
\[
\begin{aligned}
\Pr\bigl[
\deg_{\mathcal H_\Phi}(u)\ge t,\,
\deg_{\mathcal H_\Phi}(v)\ge t
\bigr]
&=
\Pr\bigl[
\deg_{\mathcal H_\Phi}(u)\ge t,\,
\deg_{\mathcal H_\Phi}(v)\ge t,\,
A
\bigr]  \\
&{}+
o\!\left(
\Pr\bigl[\deg_{\mathcal H_\Phi}(u)\ge t\bigr]
\Pr\bigl[\deg_{\mathcal H_\Phi}(v)\ge t\bigr]
\right)
\end{aligned}
\]
and 
\[
\Pr[\deg_{\mathcal H_\Phi}(u)\geq t,\deg_{\mathcal H_\Phi}(v)\geq t, A]\leq(1+o(1))\Pr[\deg_{\mathcal H_\Phi}(u)\geq t]\Pr[\deg_{\mathcal H_\Phi}(v)\geq t],
\] 
we have $\Pr[\deg_{\mathcal H_\Phi}(u)\geq t,\deg_{\mathcal H_\Phi}(v)\geq t ]\leq (1+o(1))\Pr[\deg_{\mathcal H_\Phi}(u)\geq t]\Pr[\deg_{\mathcal H_\Phi}(v)\geq t]$.
Therefore,
$\mathbb E [X^2]
\le
\mathbb E [X]+(1+o(1))(\mathbb E [X])^2.$
By \(\mathbb E X\to\infty\) and Chebyshev's inequality,
\[
\Pr[X=0]
\le
\frac{\operatorname{Var}(X)}{(\mathbb E [X])^2}\leq\frac{\mathbb{E}[X]+o((\mathbb E [X])^2)}{(\mathbb E [X])^2}
=o(1).
\]
Thus \(X>0\) with high probability, which further implies $\Delta_1(\mathcal H_\Phi)>\lambda d$ with high probability. 

\end{proof}

\begin{recalllem}{3.1,3.2,3.3}
For each iteration step $i$ in Algorithm~\ref{alg:hypergraph_container}, the following properties hold:
\begin{itemize}
    \item There always exist desirable $s_i$ and $L_i$.
    \item $\mathcal{H}_i$ is an antichain. 
    \item $\langle\mathcal{H}_i\rangle\subsetneq \langle\mathcal{H}_{i+1}\rangle$, and $\mathcal{F}_i\cap\langle\mathcal{H}_{i}\rangle=\emptyset$.
\end{itemize}  
\end{recalllem}
\begin{proof}
   For each step \(i\), we first prove the existence of a desirable pair \(s_i,L_i\).  The weight $w_p(\mathcal{H}_i^{>1})$ can be written as 
   \[w_p(\mathcal{H}_i^{>1})=\sum_{s=2}^{k}w_p(\mathcal{H}_i^s)=\sum_{s=2}^{k}\frac{p}{s}\left(\sum_{x\in C_i}w_p(\partial_x\mathcal{H}_i^s)\right).\]
   Now assume that $w_p(\mathcal{H}_i^{>1})>\delta p|C_i|$ in step $i.$
   It shows that 
   \[\delta p|C_i|< w_p(\mathcal{H}_i^{>1})=\sum_{s=2}^{k}\frac{p}{s}\left(\sum_{x\in C_i}w_p(\partial_x\mathcal{H}_i^s)\right)\leq \frac{p}{2}\left(\sum_{x\in C_i}\sum_{s=2}^{k}w_p(\partial_x\mathcal{H}_i^s)\right).\]
   there must exist $x\in C_i$ and $s\in\{2,\ldots,k\}$ such that $w_p(\partial_x\mathcal{H}_i^s)> \delta /k.$
   This guarantees the existence of $s_i$ and $L_i$. 
   
   Next, we prove by induction that each \(\mathcal H_i\) is an antichain.
The base case $i=0$ trivially holds since $\mathcal{H}_0\subseteq \mathcal{H}_{\varphi}$ and $\mathcal{H}_{\varphi}$ is an antichain.
Suppose that $\mathcal{H}_i$ is an antichain. 
The algorithm implies that for each element $F\in \mathcal{F}_i$, there must exist $F'\in \mathcal{H}_i$ such that $F\subseteq F'.$
So, no matter whether $L_i\subseteq I$ or not, $\mathcal{F}_i$ is an antichain. 
It remains to show there does not exist $E\in \mathcal{H}_i\backslash \langle\mathcal{F}_i\rangle$ and $F\in \mathcal{F}_i$ such that $E\subseteq F$ or $F\subseteq E.$ 
If such $E,F$ exist with $F\subseteq E$, then $E\in\langle\mathcal{F}_i\rangle$, a contradiction.
On the other hand, if $E, F$ exist with $E\subseteq F$, then $E\subseteq F\subsetneq F'$, contradicting the induction hypothesis that \(\mathcal H_i\) is an antichain.

It remains to prove $\mathcal{F}_i\cap\langle\mathcal{H}_{i}\rangle=\emptyset$ and $\langle \mathcal{H}_i \rangle\subsetneq \langle \mathcal{H}_{i+1}\rangle$. 
Suppose that there exists $F\in\mathcal{F}_i\cap \langle\mathcal{H}_{i}\rangle$.
Then there exists $E\in \mathcal{H}_i$ such that $E\subseteq F.$
It suffices to prove that there exists another edge $E'\in \mathcal{H}_i$ with $F\subsetneq E'$. 
Recall that in each iteration step $i$, the family $\mathcal{F}_i$ is defined as follows:
If $L_i \subseteq I$, then $\mathcal{F}_i = \partial_{L_i}\mathcal{H}_i^{s_i}$;
Otherwise, $\mathcal{F}_i = \{L_i\}$.
For the case $L_i\subseteq I,$ each $F\in \partial_{L_i}\mathcal{H}_i^{s_i}$ satisfies $L_i\cup F\in \mathcal{H}_i$.
This implies that $E'=L_i\cup F$ is our desired edge. 
For the case $L_i\not \subseteq I,$
the fact that $L_i$ is chosen implies that $\partial_{L_i}\mathcal{H}_i^{s_i}$ is not empty.
Hence, there exists $E'\in \mathcal{H}_i^{s_i}$ such that $L_i\subsetneq E'$.
Consequently, in both cases, such an edge $E'$ exists. 
This fact indicates that $E\subseteq F\subsetneq E'$, contradicts that $\mathcal{H}_i$ is an antichain.
Therefore, $\mathcal{F}_i\cap \langle\mathcal{H}_{i}\rangle=\emptyset$. 
The definition $\mathcal{H}_{i+1} = (\mathcal{H}_{i} \setminus \langle \mathcal{F}_i \rangle) \cup \mathcal{F}_i$ implies $\langle \mathcal{H}_i \rangle\subseteq \langle \mathcal{H}_{i+1}\rangle$ and $\mathcal{F}_i\subseteq \mathcal{H}_{i+1}.$
Since $\mathcal{F}_i$ is non-empty, the set $\mathcal{H}_{i+1}$ itself contains elements outside $\langle \mathcal{H}_i \rangle$.
This finally shows that $\langle \mathcal{H}_i \rangle\subsetneq \langle \mathcal{H}_{i+1}\rangle,$ completing the proof.
\end{proof}

\begin{recalllem}{3.11}
For each independent set $I\in \mathcal{I}(\mathcal{H})$, $S_J\subseteq I\subseteq C_J$, where $J$ is the stopping index in Algorithm~\ref{alg:hypergraph_container}.
Moreover, for two different independent sets $I$ and $I'$ of $\mathcal{H}_0$, if $S=S'$, then $C=C'$ and $\mathcal{G}=\mathcal{G}'$.
\end{recalllem}
\begin{proof} 
If $I$ is an independent set in $\mathcal{H}_J$, then it does not contain 1-edges, and thus $I$ is contained in $C_{J}$. 
Hence, it suffices to prove that for each iteration $i$, $I\in \mathcal{I}(\mathcal{H}_i)$ always holds.
The base case $I\in \mathcal{I}(\mathcal{H}_0)$ trivially holds. 
Assume that $I\in \mathcal{I}(\mathcal{H}_i)$. 
In the iteration process, all possible edges added to $\mathcal{H}_{i+1}$ belong to $\mathcal{F}_i$, which are not contained in $I$. 
This implies that $I$ does not contain any edge in $\mathcal{H}_{i+1}$.

It remains to compare $(C,\mathcal{G})$ and $(C',\mathcal{G}')$ if two independent sets share the same fingerprint. 
For each $i\leq J$, we will prove $S_i=S_i'$ and $\mathcal{H}_i=\mathcal{H}_i'$ by induction on $i$. 
For the base case, $S_0=S_0'=\emptyset,$ and $\mathcal{H}_0=\mathcal{H}_0=\mathcal{H}'$ hold. 
    Suppose that $S_i=S_i'$ and $\mathcal{H}_i=\mathcal{H}_i'$ hold. It implies that $C_i=C_i', s_i=s_i'$ and $L_i=L_i'$. 
    If $L_i\subseteq I$, then $S_{i+1}=S_i\cup L_i$ and $$L_i'=L_i\subseteq S_{i+1}\subseteq S=S'\subseteq I',$$ 
    where $S'\subseteq I'$ holds because $\emptyset =S'_0\subseteq I'$ and in each step we only add vertices in $I'$ to $S'$.
    This implies that $L_i'\subseteq I'$ and $S_{i+1}=S_{i+1}'=S_i\cup L_i$.
    Otherwise, both $L_i\not\subseteq I$ and $L_i'\not\subseteq I'$ (because $I$ and $I'$ are symmetric, the case $L_i'\subseteq I'$ also implies $L_i\subseteq I$).
    This implies that $S_{i+1}=S_{i+1}'=S_i.$
    Hence, in both cases, $S_{i+1}=S_{i+1}'$ and $\mathcal{F}_i=\mathcal{F}_i'$.
    This implies $\mathcal{H}_{i+1}=\mathcal{H}_{i+1}'.$
    As a consequence, $S_J=S_J'$, $C_J=C_J'$ and $\mathcal{H}_J=\mathcal{H}_J'.$
    The property that whether $C_J\cap \{x_i,\neg x_i\}=\emptyset$ is the same as $C_J',$ which shows that either  $C_J$ and $C_J'$ are both included in $\mathcal{C}$, or  are both deleted. 
    Hence, the existence of the final output $(S,C,\mathcal{G})$ is the same with $(S',C',\mathcal{G'})$.
\end{proof}

\end{document}